\documentclass[11pt]{amsart}
\usepackage{geometry}                
\usepackage[parfill]{parskip}    
\usepackage{graphicx}
\usepackage{amssymb}
\usepackage{epstopdf}
\usepackage{subfigure}
\usepackage{amsmath}
\usepackage{color}
\swapnumbers 
\numberwithin{equation}{section} 

\theoremstyle{plain} 
\newtheorem{thm}{Theorem}[section] 
\newtheorem{cor}[thm]{Corollary} 
 
\newtheorem{lem}[thm]{Lemma} 
 
\newtheorem{prop}[thm]{Proposition}

\theoremstyle{definition}

\def\Z{{\mathbb Z}} 
\def\fh{{\hat{\pi}}}
\def\sgn{\protect\operatorname{sgn}}
\def\im{\protect\operatorname{Im}}

\title[A Generalization of the Turaev Cobracket]{A Generalization of the Turaev Cobracket and the Minimal Self-Intersection Number of a Curve on a Surface}
\author{Patricia Cahn}
\date{}                                           

\begin{document}
\maketitle
\begin{abstract} 
 Goldman and Turaev constructed a Lie bialgebra structure on the free $\mathbb{Z}$-module generated by free homotopy classes of loops on a surface.  Turaev conjectured that his cobracket $\Delta(\alpha)$ is zero if and only if $\alpha$ is a power of a simple class.  Chas constructed examples that show Turaev's conjecture is, unfortunately, false.  We define an operation $\mu$ in the spirit of the Andersen-Mattes-Reshetikhin algebra of chord diagrams.  The Turaev cobracket factors through $\mu$, so we can view $\mu$ as a generalization of $\Delta$.  We show that Turaev's conjecture holds when $\Delta$ is replaced with $\mu$.  We also show that $\mu(\alpha)$ gives an explicit formula for the minimum number of self-intersection points of a loop in $\alpha$.  The operation $\mu$ also satisfies identities similar to the co-Jacobi and coskew symmetry identities, so while $\mu$ is not a cobracket, $\mu$ behaves like a Lie cobracket for the Andersen-Mattes-Reshetikhin Poisson algebra.
\end{abstract} 

\section{Introduction}
\noindent We work in the smooth category.  All manifolds and maps are assumed to be smooth unless stated otherwise, where smooth means $C^{\infty}$.

Goldman \cite{Goldman} and Turaev \cite{Turaev} constructed a Lie bialgebra structure on the free $\Z$-module generated by nontrivial free homotopy classes of loops on a surface $F$.  Turaev \cite{Turaev} conjectured that his cobracket $\Delta(\alpha)$ is zero if and only if the class $\alpha$ is a power of a simple class, where we say a free homotopy class is simple if it contains a simple representative.  Chas \cite{Chas} constructed examples showing that, unfortunately, Turaev's conjecture is false on every surface of positive genus with boundary.  In this paper, we show that Turaev's conjecture is almost true.   We define an operation $\mu$ in the spirit of the Andersen-Mattes-Reshetikhin algebra of chord diagrams, and show that Turaev's conjecture holds on all surfaces when one replaces $\Delta$ with $\mu$.  
\begin{figure}[htbp] 
   \centering
   \includegraphics[width=4in]{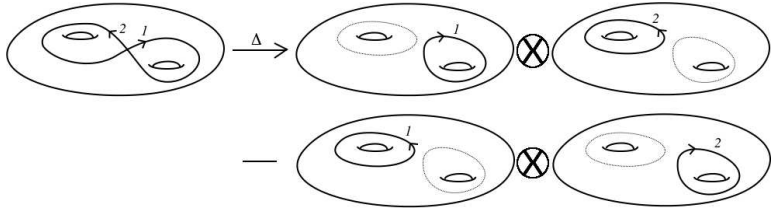} 
   \caption{Two terms of Turaev's cobracket $\Delta(\alpha)$ with coefficients $+1$ and $-1$.}
   \label{newintrocobracket.fig}
\end{figure} 

Turaev's cobracket $\Delta(\alpha)$ is a sum over the self-intersection points $p$ of a loop $a$ in a free homotopy class $\alpha$.  Each term of the sum is a simple tensor of free homotopy classes loops, which are obtained by smoothing $a$ at the self-intersection point $p$ along its orientation.  Each simple tensor is equipped with a sign coming from the intersection at $p$ (see Figure \ref{newintrocobracket.fig}).  Turaev's conjecture is false because it is not uncommon for the same simple tensor of loops to appear twice in the sum $\Delta(\alpha)$, but with different signs.  

\begin{figure}[htbp] 
   \centering
   \includegraphics[width=4in]{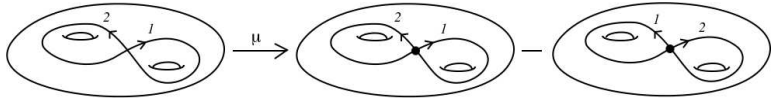} 
   \caption{Two terms of the operation $\mu(\alpha)$ with coefficients $+1$ and $-1$.}
\end{figure}
\begin{figure}[h]\center
\scalebox{1.0}{\includegraphics[width=4cm]{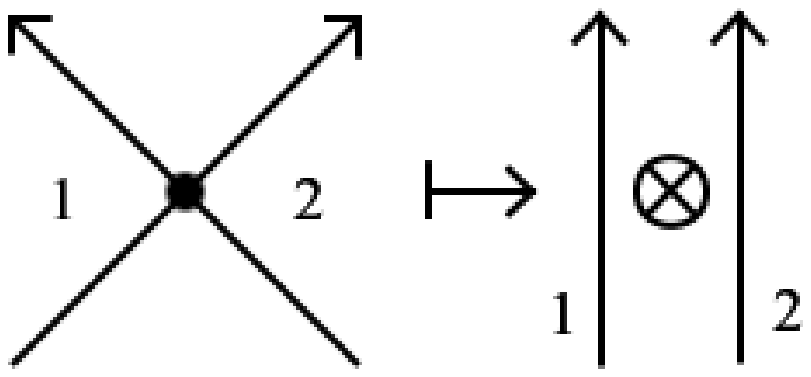}}  
\label{factors1.fig}
\end{figure}

We define the operation $\mu(\alpha)$ as a sum over the self-intersection points $p$ of a loop $a$ in $\alpha$, as in the definition of the Turaev cobracket.  Rather than smoothing at each self-intersection point to obtain a simple tensor of two loops, we glue those loops together to create a wedge of two circles mapped to the surface.  This can also be viewed as a chord diagram with one chord.  As a result, terms of $\mu$ are less likely to cancel than terms of $\Delta$, and hence $\mu(\alpha)$ is less likely to be zero.  In fact, Turaev's conjecture holds when formulated for $\mu$ rather than $\Delta$:
\begin{thm}\label{main2} Let $F$ be an oriented surface with or without boundary, which may or may not be compact.  Let $\alpha$ be a free homotopy class on $F$.  Then $\mu(\alpha)=0$ if and only if $\alpha$ is a power of a simple class.  
\end{thm}
There is a simple relationship between $\Delta$ and $\mu$; namely, if one smoothes each term of $\mu$ at the gluing point, and tensors the resulting loops, one obtains a term of $\Delta$ (see Figure \ref{factors1.fig}).
  Hence the Turaev cobracket factors through $\mu$, and we can view $\mu$ as a generalization of $\Delta$.  The relationship between $\mu$ and $\Delta$ is analogous to the relationship between the Andersen-Mattes-Reshetikhin Poisson bracket for chord diagrams and the Goldman Lie bracket.  It is natural to wonder to what extent we can view $\mu$ as a cobracket for the Andersen-Mattes-Reshetikhin algebra.  While $\mu$ is not a cobracket, in the final section of the paper, we show that $\mu$ satisfies identities similar to coskew symmetry and the co-Jacobi identity.\\\\
The operation $\mu$ also gives an explicit formula for the minimum number of self-intersection points of a generic loop in a given free homotopy class $\alpha$.  We call this number the {\it minimal self-intersection number of $\alpha$} and denote it by $m(\alpha)$.  Both Turaev's cobracket and the operation $\mu$ give lower bounds on the minimal self-intersection number of a given homotopy class $\alpha$.  We call a free homotopy class {\it primitive} if it is not a power of another class in $\pi_1(F)$.  Any class $\alpha$ can be written as $\beta^n$ for some primitive class $\beta$ and $n\geq 1$.  It follows easily from the definitions of $\Delta$ and $\mu$ that $m(\alpha)$ is greater than or equal to $n-1 $ plus half the number of terms in the (reduced) linear combinations $\Delta(\alpha)$ or $\mu(\alpha)$.  More formally, the number of terms $t(L)$ of a reduced linear combination $L$ of simple tensors of classes of loops, or of classes chord diagrams, is the sum of the absolute values of the coefficients of the classes.  Chas' counterexamples to Turaev's conjecture show that the lower bound given by $\Delta(\alpha)$ cannot, in general, be used to compute the minimal self-intersection number of $\alpha$.  However the lower bound given by $\mu(\alpha)$ is always equal to $m(\alpha)$:
\begin{thm}\label{nonprimthm}  Let $F$ be an oriented surface with or without boundary, which may or may not be compact.  Let $\alpha$ be a nontrivial free homotopy class on $F$ such that $\alpha=\beta^n$, where $\beta$ is primitive and $n\geq 1$.  Then the minimal self-intersection number of $\alpha$ is equal to $n-1$ plus the half number of terms of $\mu(\alpha)$.
\end{thm}
In order to prove the case of Theorem \ref{nonprimthm} where $n>1$, we make use of the results of Hass and Scott \cite{HassScott} who describe geometric properties of curves with minimal self-intersection (see also \cite{FreedmanHassScott}).  

We briefly summarize some results related to Turaev's conjecture and to computing the minimal self-intersection number.  Le Donne \cite{LeDonne} proved that Turaev's conjecture is true for genus zero surfaces.  For surfaces of positive genus, one might wonder to what extent Turaev's conjecture is false.  Chas and Krongold \cite{ChasKrongold} approach this question by showing that, on surfaces with boundary, if $\Delta(\alpha)=0$ and $\alpha$ is at least a third power of a primitive class $\beta$, then $\beta$ is simple. \\\\
A nice history of the problem of determining when a homotopy class is represented by a simple loop is given in Rivin \cite{Rivin}.  Birman and Series \cite{BirmanSeries} give an explicit algorithm for detecting simple classes on surfaces with boundary.  Cohen and Lustig \cite{CohenLustig} extend the work of Birman and Series to obtain an algorithm for computing the minimal intersection and self-intersection numbers of curves on surfaces with boundary, and Lustig \cite{Lustig} extends this to closed surfaces.  We give an example which shows how one can algorithmically compute $m(\alpha)$ using $\mu$ on surfaces with boundary, though generally we do not emphasize algorithmic implications in this paper. \\\\
A different algebraic solution to the problem of computing the minimal intersection and self-intersection numbers of curves on a surface is given by Turaev and Viro \cite{TuraevViro}.  The advantage of $\mu$ is that it has a simple relationship to $\Delta$ and pairs well with the Andersen-Mattes-Reshetikhin Poisson bracket.  In fact, Chernov \cite{Chernov} uses the Andersen-Mattes-Reshetikhin bracket to compute the minimum number of intersection points of loops in given free homotopy classes.

\section{The Goldman-Turaev and Andersen-Mattes-Reshetikhin Algebras and the Operation $\mu$}
\subsection{The Goldman-Turaev Lie Bialgebra}  We will now define the Goldman-Turaev Lie Bialgebra on the free $\Z$-module generated by the set $\fh$ of free homotopy classes of loops on $F$, which we denote by $\Z[\fh]$.  Let $\alpha, \beta \in \fh$, and let $a$ and $b$ be smooth, transverse representatives of $\alpha$ and $\beta$, respectively.  We will use square brackets to denote the free homotopy class of a loop.  The set of intersection points $\mathcal{I}_{a, b}$, or just $\mathcal{I}$ when the choice of $a$ and $b$ is clear, is defined to be
$$\mathcal{I}=\{(t_1,t_2)\in S^1\times S^1: a(t_1)=b(t_2)\}.$$
Let $a\cdot_p b$ denote the product of $a$ and $b$ as based loops in $\pi_1(F,p)$, where $p=a(t_1)=b(t_2)$ and $(t_1,t_2)\in \mathcal{I}$.  If $p$ is the image of more than one ordered pair in $\mathcal{I}$, then there is more than one homotopy class in $\pi_1(F,p)$ corresponding to $a$ (or to $b$), so we choose a class as follows:  Let $a_*:\pi_1(S^1,t_1)\rightarrow\pi_1(F, a(t_1))$ be the induced map on the fundamental groups.  Let $\gamma=[g(t)]$ be the generator of $\pi_1(S^1,t_1)$ whose orientation agrees with the chosen orientation of $S^1$.  Then the class of $a$ in $\pi_1(F,p)$ is given by $a_*[g(t)]$.  We choose the class of $b$ in $\pi_1(F,p)$ in the same way.  In particular, we must specify preimages of $p$ under $a$ and $b$ (i.e. a point in $\mathcal{I}$) for the notation $a\cdot_p b$ to make sense.\\\\
The Goldman bracket \cite{Goldman} is a linear map $[\cdot,\cdot]:\Z[\fh]\otimes_{\Z}\Z[\fh]\rightarrow \Z[\fh]$, defined by
$$[\alpha,\beta]=\sum_{(t_1,t_2)\in \mathcal{I} } \sgn(p;a,b)[a\cdot_p b],$$
where $\sgn(p;a,b)=1$ if the orientation given by the pair of vectors $\{\bf{a'}(t_1), \bf{b'}(t_2)\}$ agrees with the orientation of $F$, and $\sgn(p;a,b)=-1$ otherwise.  To check that the definition of $[\cdot,\cdot]$ is independent of the choices of $a$ and $b$, one must show that $[\alpha,\beta]$ does not change under elementary moves for a pair of smooth curves in general position.  Using linearity, the definition of $[\cdot,\cdot]$ can be extended to all of $\Z[\fh]\otimes_{\Z}\Z[\fh]$.\\\\
Next we define the Turaev cobracket \cite{Turaev}.  Let $\alpha$ be a free homotopy class on $F$, and let $a$ be a smooth representative of $\alpha$ with transverse self-intersection points.  Let $\mathcal{SI}_{a}$, or just $\mathcal{SI}$ when the choice of $a$ is clear, denote the set of self-intersection points of the loop $a$.  Let $\mathbb{D}$ be the diagonal in $S^1\times S^1$.  Elements of $\mathcal{SI}$ will be points in $S^1\times S^1-\mathbb{D}$ modulo the action of $\Z_2$ which interchanges the two coordinates.  Now we define
$$\mathcal{SI}=\{(t_1,t_2)\in (S^1\times S^1-\mathbb{D})/\Z_2:a(t_1)=a(t_2)\}.$$
Let $p=a(t_1)=a(t_2)$ be a self-intersection point of $a$.  Let $[t_1,t_2]$ denote the arc of $S^1$ going form $t_1$ to $t_2$ in the direction of the orientation of $S^1$, and let $[t_2,t_1]$ denote the arc of $S^1$ going from $t_2$ to $t_1$ in the direction of the orientation of $S^1$.  Since $p=a(t_1)=a(t_2)$, then $a([t_1,t_2])$ and $a([t_2,t_1])$ are loops.  We assign these loops the names $a^1_p$ and $a^2_p$ in such a way that the ordered pair of tangent vectors $\{(a^1_p)',(a^2_p)'\}$ gives the chosen orientation of $T_pF$.  Now we let $\mathcal{SI}_0$ be the subset of $\mathcal{SI}$ which contains only self-intersection points $p$ such that the loops $a^i_p$ are nontrivial:
$$\mathcal{SI}_0=\{ (t_1,t_2) \in\mathcal{SI}:p=a(t_1)=a(t_2), a^1_p , a^2_p \neq 1\in \pi_1(F_p)\}.$$
The Turaev cobracket is a linear map $\Delta:\Z[\fh]\rightarrow \Z[\fh]\otimes_{\Z} \Z[\fh]$ which is given on a single homotopy class by
$$\Delta(\alpha) = \sum_{(t_1,t_2)\in \mathcal{SI}_0}[a^1_p]\otimes[a^2_p]-[a^2_p]\otimes[a^1_p].$$
One can show that the definition of $\Delta$ is independent of the choice of $a\in \alpha$ by showing $\Delta(\alpha)$ does not change under elementary moves for a smooth loop in general position.  Using linearity, this definition of $\Delta$ can be extended to all of $\Z[\fh]$.\\\\
Together, $[\cdot,\cdot]$ and $\Delta$ equip $\Z[\fh]$ with an involutive Lie Bialgebra structure \cite{Goldman, Turaev}.  That is, $[\cdot,\cdot]$ and $\Delta$ satisfy (co)skew-symmetry, the (co) Jacobi identity, a compatibility condition, and $[\cdot,\cdot]\circ\Delta=0$.  A complete definition of a Lie Bialgebra is given in \cite{Chas}.

\subsection{The Andersen-Mattes-Reshetikhin Algebra of Chord Diagrams}  
\begin{figure}[h]\center
\scalebox{1.0}{\includegraphics[width=8cm]{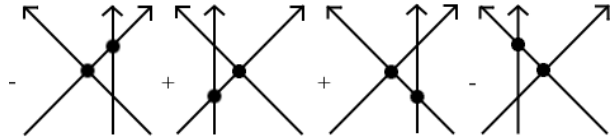}}  
\caption{$4T$-relations}\label{4term.fig}
\end{figure}

\begin{figure}[h]\center
\scalebox{1.0}{\includegraphics[width=4cm]{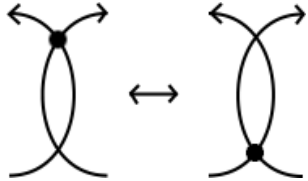}}  
\caption{An elementary move for chord diagrams (with one of several possible choices of orientations on the arcs).}\label{type2movewithdot.fig}
\end{figure}

\begin{figure}[h]\center
\scalebox{1.0}{\includegraphics[width=5cm]{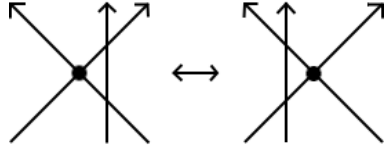}}  
\caption{An elementary move for chord diagrams (with one of several possible choices of orientations on the arcs).}\label{type3movewithdot.fig}
\end{figure}

We now summarize the Andersen-Mattes-Reshetikhin algebra of chord diagrams on $F$ \cite{AMR2, AMR}.  A \textit{chord diagram} is a disjoint union of oriented circles $S_1,...,S_k$, called \textit{core circles}, along with a collection of disjoint arcs $C_1,...,C_l$, called \textit{chords},  such that\\
1) $\partial C_i \bigcap \partial C_j=\emptyset$ for $i\neq j$, and\\
2) $\bigcup_{i=1}^{l} \partial C_i = \left(\bigcup_{i=1}^{k}S_i\right) \bigcap \left( \bigcup_{i=1}^{l}C_i\right)$.\\\\\
A \textit{geometrical chord diagram on $F$} is a smooth map from a chord diagram $D$ to $F$ such that each chord $C_i$ in $D$ is mapped to a point.  A \textit{chord diagram on $F$} is a homotopy class of a geometrical chord diagram $D$, denoted $[D]$.\\\\
Let $M$ denote the free $\Z$-module generated by the set of chord diagrams on $F$ (\cite{AMR} uses coefficients in $\mathbb{C}$, but we use $\Z$ here for consistency).  Let $N$ be the submodule generated by a set of $4T$-relations, one of which is shown in Figure \ref{4term.fig}.  The other relations can be obtained from this one as follows: one can reverse the direction of any arrow, and any time a chord intersects an arc whose orientation is reversed, the diagram is multiplied by a factor of -1.\\
 \\\\
Given two chord diagrams $\mathcal{D}_1$ and $\mathcal{D}_2$ on $F$, we can form their disjoint union by choosing representatives (i.e., geometrical chord diagrams) $D_i$ of $\mathcal{D}_i$, taking a disjoint union of their underlying chord diagrams, mapping the result to $F$ as prescribed by the $D_i$, and taking its free homotopy class. The disjoint union of chord diagrams $\mathcal{D}_1 \cup \mathcal{D}_2$ defines a commutative multiplication on $M$, giving $M$ an algebra structure with $N$ as an ideal.  Let $ch=M/N$, and call this the \textit{algebra of chord diagrams}.  \\\\
Andersen, Mattes, and Reshetikhin \cite{AMR2, AMR} constructed a Poisson bracket on $ch$, which can be viewed as a generalization of the Goldman bracket for chord diagrams on $F$ rather than free homotopy classes of loops.  Let $\mathcal{D}_1$ and $\mathcal{D}_2$ be chord diagrams on $F$, and choose representatives $D_i$ of $\mathcal{D}_i$.  We define the set of intersection points $\mathcal{I}_{D_1,D_2}$, or just $\mathcal{I}$ when the choice of $D_1$ and $D_2$ is clear, to be
$$\mathcal{I}=\{ (t_1,t_2): D_1(t_1)=D_2(t_2)\}, $$
where $t_i$ is a point in the preimage of the geometrical chord diagrams $D_i$.  For each $(t_1,t_2)\in\mathcal{I}$ with $p=D_i(t_i)$, let $D_1 \cup_p D_2$ denote the geometrical chord diagram obtained by adding a chord between $t_1$ and $t_2$.  It is necessary to specify preimages of $p$ for this notation to be well-defined.  Since each copy of $S^1$ in the chord diagram is oriented, we can define $\sgn(p; D_1,D_2)$ as before.  The Andersen-Mattes-Reshetikhin Poisson bracket $\{\cdot,\cdot\}:ch\times ch \rightarrow ch$ is defined by
$$\{\mathcal{D}_1,\mathcal{D}_2\}=\sum_{(t_1,t_2)\in \mathcal{I}}\sgn(p;D_1,D_2) [D_1 \cup_p D_2],$$
where square brackets denote the free homotopy class of a geometrical chord diagram.  This definition of $\{\cdot,\cdot\}$ can be extended to all of $ch$ using bilinearity.  For a proof that $\{\cdot,\cdot\}$ does not depend on the choices of $D_i\in\mathcal{D}_i$, $i=1,2$, see \cite{AMR}.  In particular, it is necessary to check that $\{\cdot,\cdot\}$ is invariant under elementary moves, including the Reidemeister moves and the moves in Figures \ref{type2movewithdot.fig} and \ref{type3movewithdot.fig}, and the $4T$-relations.

\subsection{The Operation $\mu$}  
\begin{figure}[b]\center
\scalebox{1.0}{\includegraphics[width=4cm]{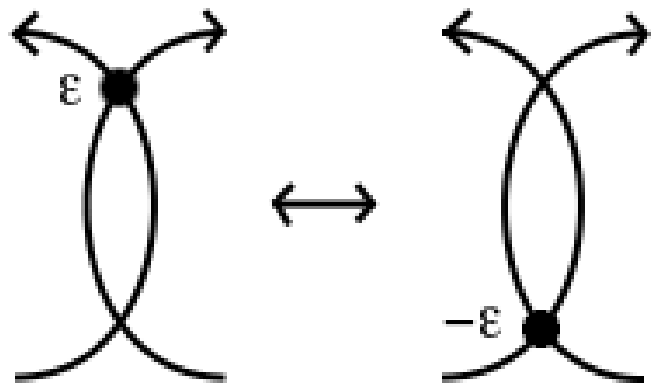}}  
\caption{Elementary move for chord diagrams with oriented chords.}\label{type2movewithdotandsign.fig}
\end{figure}

\begin{figure}[b]\center
\scalebox{1.0}{\includegraphics[width=5cm]{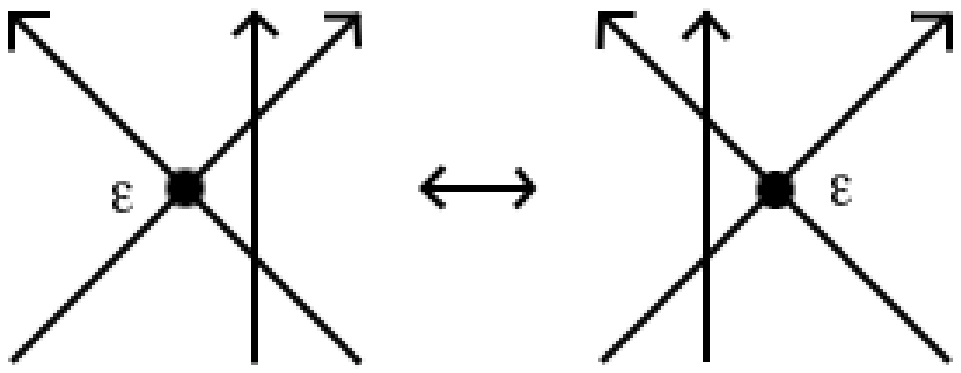}}  
\caption{Elementary move for chord diagrams with oriented chords.}\label{type3movewithdotandsign.fig}
\end{figure}

The definition of $\mu$ given in this section is the simplest for the purposes of computing the minimal self-intersection number of a free homotopy class $\alpha$.  In this section, we define $\mu$ only on free homotopy classes.  In the final section of this paper, we modify the definition of $\mu$ in a way that allows us to more easily state an analogue of the co-Jacobi identity, and which allows us to extend the definition of $\mu$ to certain chord diagrams in the Andersen-Mattes-Reshetikhin algebra.  The modified definition agrees with the definition below for free homotopy classes.\\\\
 For this defintion of $\mu$, we will need to use chord diagrams with oriented chords.  Suppose $C$ is an oriented chord with its tail at $t\in S^1$ and its head at $h\in S^1$ in a geometrical chord diagram $D$.  We say $C$ \textit{agrees with the orientation of $F$} if the ordered pair of vectors $\{{\bf D'}(h), {\bf D'}(t)\}$ gives the chosen orientation of $F$.  When we draw the image of a geometrical chord diagram, we label the image of a chord $C$ with a $`+$' if $C$ agrees with the orientation of $F$, and we label it with a `$-$' otherwise.\\\\
Let $E$ denote the free $\Z$-module generated by chord diagrams on $F$ consisting of one copy of $S^1$, and one oriented chord connecting distinct points of that copy of $S^1$.  In addition to the usual Reidemeister moves, we have two additional elementary moves for diagrams with signed chords.   These moves, with one possible choice of orientation on the branches, are shown in Figures \ref{type2movewithdotandsign.fig} and \ref{type3movewithdotandsign.fig}, where $\epsilon \in \{+,-\}$ denotes the sign on the chord.  We define a linear map $\mu: \Z[\fh]\rightarrow E$.  \\\\
Let $D:S^1\rightarrow F$ be a geometrical chord diagram on $F$ with one core circle.  For each self-intersection point $p=D(t_1)=D(t_2)$ of $D$, we let $D^+_p$ (respectively $D^-_p$) be the geometrical chord diagram obtained by adding an oriented chord between $t_1$ and $t_2$ that agrees (respectively, does not agree) with the orientation of $F$.  \\\\
Now we define $\mu$ on the class of the geometrical chord diagram $D$ by
$$\mu([D])=\sum_{(t_1,t_2)\in \mathcal{SI}_0}[D^+_p]-[D^-_p].$$
Using linearity, we can extend this definition to all of $\Z[\fh]$.  It remains to check that $\mu([D])$ is independent of the choice of representative of $[D]$.\\

\subsection{$\mu(\mathcal{D})$ is independent of the choice of representative of $\mathcal{D}$} We check that $\mu$ is invariant under the usual Reidemeister moves:
\begin{enumerate}
\item {\it Regular isotopy:} Invariance is clear.
\item {\it First Reidemeister Move:} This follows from the definition of $\mathcal{SI}_0$.
\item {\it Second Reidemeister Move:} This follows from the move in Figure \ref{type2movewithdotandsign.fig}.
\item{\it Third Reidemeister Move:} This follows from the move in Figure \ref{type3movewithdotandsign.fig}.
\end{enumerate}
We note that when checking invariance under the second and third moves, one must consider the case where some of the self-intersection points are in $\mathcal{SI}$ but not in $\mathcal{SI}_0$.

\subsection{Alternative notation for $\mu$}  We would like to show that $\Delta$ factors through $\mu$.  To do this, we will rewrite the definition of $\mu$ for a free homotopy class in a way that makes its relationship to $\Delta$ more transparent.  Let $\phi$ and $\psi: I=[0,1]\rightarrow F$ be loops in $F$ based at $p$, such that $ \phi'(0)=\psi'(1)$ and $\phi'(1)=\psi'(0)$.  We define a geometrical chord diagram $\phi \bullet_p \psi$ which glues the loops $\phi$ and $\psi$ at the point $p$.   The underlying chord diagram of $\phi \bullet_p \psi$ contains one core circle $S^1=I/\partial I$, and one oriented chord $C$ with its head at $0\in I$ and its tail at $\frac12\in I$.  The geometrical chord diagram $\phi \bullet_p \psi$ maps the chord $C$ to $p$.  Then we define $(\phi \bullet_p \psi) |_{[0,\frac12]}=\phi$ and $(\phi \bullet_p \psi) |_{[\frac12,1]}=\psi$.  \\\\
Now we are ready to rewrite the definition of $\mu$ for $\alpha\in \fh(F)$.  Let $a$ be a representative of $\alpha$, and for each $(t_1,t_2)\in\mathcal{SI}_0$ with $p=a(t_1)=a(t_2)$, let $a^1_p$ and $a^2_p$ be the loops we defined for the Turaev cobracket.  Now
$$\mu(\alpha)=\sum_{(t_1,t_2)\in \mathcal{SI}_0}[a^1_p\bullet_p a^2_p]-[a^2_p\bullet_p a^1_p].$$

\subsection{Relationship between $\mu$, the Goldman-Turaev Lie bialgebra, and the Andersen-Mattes-Reshetikhin Algebra of Chord Diagrams}
\begin{figure}[h]\center
\scalebox{1.0}{\includegraphics[width=4cm]{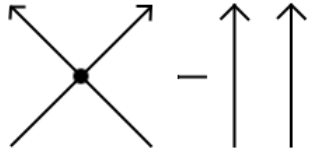}}  
\caption{Generator of $I$}\label{I.fig}
\end{figure}

\begin{figure}[h]\center
\scalebox{1.0}{\includegraphics[width=4cm]{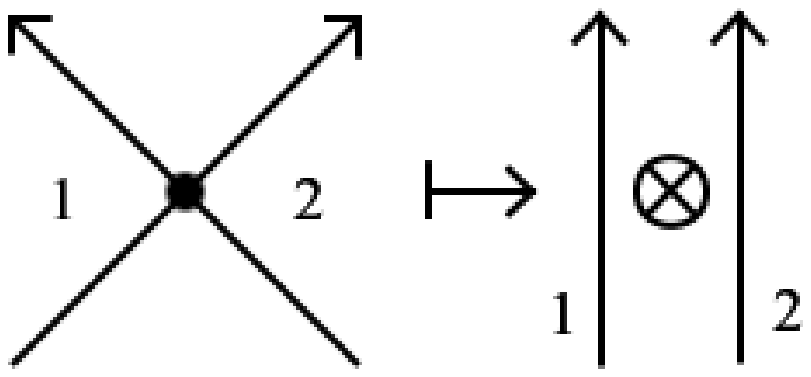}}  
\caption{The map $Q$}\label{factors.fig}
\end{figure}

Andersen, Mattes and Reshetikhin \cite{AMR} show that there is a quotient algebra of $ch$ which corresponds to Goldman's algebra.  Let $I$ be the ideal generated by the relation in Figure \ref{I.fig}.  In the quotient $ch/I$, each chord diagram is identified with the disjoint union of free homotopy classes obtained by smoothing the diagram at the intersections which are images of chords. 
One can check that $P: ch\rightarrow ch/I$ is a Poisson algebra homomorphism and $ch/I$ is a Poisson algebra with an underlying Lie algebra that corresponds to Goldman's algebra \cite{AMR}. \\\\
There is a similar relationship between the Turaev cobracket and $\mu$.  Let $Q$ be the map which smoothes the chord diagram according to its orientation at an intersection which is an image of a chord, and tensors the two resulting homotopy classes together (see Figure \ref{factors.fig}). Then $\Delta=Q\circ \mu$.\\\\
\textit{Remark:}  Turaev \cite[p. 660]{Turaev} notes that the Turaev cobracket can be obtained algebraically from an operation defined in Supplement 2 of \cite{TuraevIntersection}.  It is possible that $\mu$ may be obtained from this operation as well.  We do not know a way of obtaining Turaev's operation from $\mu$.

\section{Proofs of Theorems}
\noindent In this section, we prove Theorems \ref{main2} and \ref{nonprimthm}.  Recall that Theorem \ref{main2} states that $\mu(\alpha)=0$ if and only if $\alpha$ is a power of a simple class.  Theorem \ref{nonprimthm} gives an explicit formula for $m(\alpha)$.  We begin by describing two types of self-intersection points of a loop which is freely homotopic to a power of another loop.  Then we prove Theorem \ref{mainlemma}, which describes when certain terms of $\mu$ cancel.  Theorems \ref{main2} and \ref{nonprimthm} are corollaries of Theorem \ref{mainlemma}. 
\subsection{Intersection Points of Powers of Loops}  
\begin{figure}[h]\center
\scalebox{1.0}{\includegraphics[width=8cm]{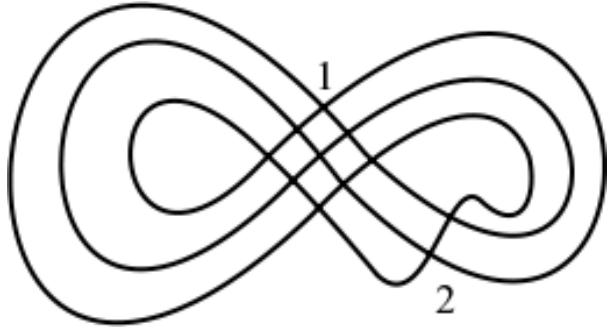}}  
\caption{Type 1 and Type 2 self-intersection points.}\label{type1type2intpts.fig}
\end{figure}

\begin{figure}[h]\center
\scalebox{1.0}{\includegraphics[width=8cm]{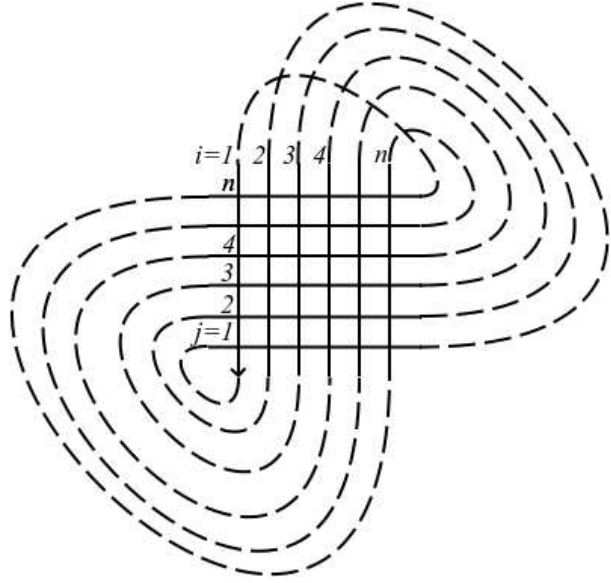}}  
\caption{Type 1 intersection points.}\label{powerslabelled.fig}
\end{figure}

Our goal is to understand the conditions under which different terms of $\mu(\alpha)$ cancel, when $\alpha\in \fh$ is a power of another class $\beta$ in $\pi_1(F)$.  To do this, we need to distinguish between two different types of self-intersection points of a curve.  Suppose we choose a geodesic representative $g$ of $\alpha$.  Either all self-intersection points of $g$ are transverse, or $g$ has infinitely many self-intersection points, and in particular, $g$ is a power of another geodesic.  Let $p$ be a point on the image of $g$ which is not a transverse self-intersection point of $g$.  Let $h$ be a geodesic loop such that $g=h^n$ in $\pi_1(F,p)$, and such that there is no geodesic $f$ such that $h=f^k$ (it is possible that $|n|=1$).  Now we know that $h$ has finitely many self-intersection points, all of which are transverse.  Let $m$ be the number of self-intersection points of $h$.  Since $F$ is orientable, we can perturb $g$ slightly to obtain a loop $g'$ as follows: We begin to traverse $g$ beginning at $p$, but whenever we are about to return to $p$, we shift slightly to the left.  After doing this $n$ times, we must return to $p$ and connect to the starting point.  This requires crossing $n-1$ strands of the loop, creating $n-1$ self-intersection points.   We call these \textit{Type 2 self-intersection points}.  For self-intersection point of $h$, we get $n^2$ self-intersection points of $g$ (see Figure \ref{type1type2intpts.fig}).  We call these $mn^2$ self-intersection points \textit{Type 1 self-intersection points}.  We note that we are counting self-intersections with multiplicity, as some of the self-intersection points of $h$ may be images of multiple points in $\mathcal{SI}$.
Given a transverse self-intersection point $p$ of $h$, we will denote the corresponding set of $n^2$ Type 1 self-intersection points of $g'$ by $\{p_{i,j}\}$, where $i\in\{1,...,n\}$ is the label on the strand corresponding first branch of $h$ at $p$ (i.e., a strand going from top to bottom in Figure \ref{powerslabelled.fig}), and $j\in\{1,...,n\}$ is the label on the strand corresponding to the second branch of $h$ at $p$ (i.e., a strand going from left to right in Figure \ref{powerslabelled.fig}).  This relationship between the numbers of self-intersection points of $g$ and $h$ can be found in \cite{TuraevViro} for both orientable and non-orientable surfaces.
\begin{lem} \label{type1lemma} Let $g$ be a geodesic representative of $\alpha \in \fh(F)$, with $g=h^n$, and $h$ and $n$ are as defined in the paragraph above.  Then the contribution to $\mu(\alpha)$ of a Type 1 self-intersection point $p_{i,j}$ is
 $$[(XY)^I X\bullet_p (YX)^JY]-[(YX)^JY\bullet_p (XY)^IX], $$
 where $X=h^1_p$, $Y=h^2_p$, and $I,J\in \mathbb{N}$ such that $I+J=n-1$.
\end{lem}
\noindent\textit{Proof.}  We will compute the contribution to $\mu$ for a Type 1 self-intersection point $p_{i,j}$ of $g'$, where $g'$ is the perturbed version of $g$ described in the above paragraph.  These terms are $[(g')^1_{p_{i,j}}\bullet_{p_{i,j}} (g')^2_{p_{i,j}}]$ and $-[(g')^2_{p_{i,j}}\bullet_{p_{i,j}} (g')^1_{p_{i,j}}]$.  However, when we record the terms of $\mu$, we perturb $g'$ back to $g$, so that the terms we record are geometrical chord diagrams whose images are contained in the image of $g$ and whose chords are mapped to $p$.  To compute $(g')^1_{p_{i,j}}$, we begin at $p_{i,j}$ along the branch corresponding to $X=h^1_p$, and wish to know how many times we traverse branches corresponding to $X=h^1_p$ and $Y=h^2_p$ before returning to $p_{i,j}$.  The first time we return to $p_{i,j}$, we must return along the $j^{th}$ branch of $X$.  Therefore $[(g')^1_{p_{i,j}}]=[(XY)^IX]$ for some integer $I\geq 0$.  If we begin at $p_{i,j}$ along the branch corresponding to $Y=h^2_p$, we return to $p_{i,j}$ for the first time on the $i^{th}$ branch of $Y$.  Therefore $[(g')^2_{p_{i,j}}]=[(YX)^JY]$ for some integer $J\geq 0$.  But if we traverse $(g')^1_{p_{i,j}}$ followed by $(g')^2_{p_{i,j}}$, we must traverse $g'$ exactly once, so $I+J=n-1$.
\qed

\subsection{Canceling terms of $\mu$}  
\begin{figure}[h]\center
\scalebox{1.0}{\includegraphics[width=10cm]{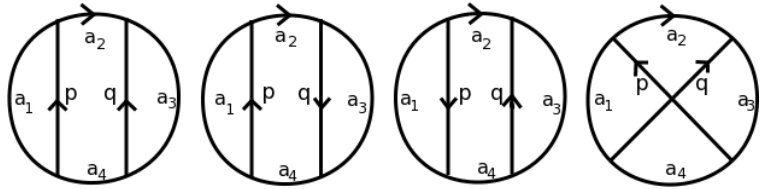}}  
\caption{$(a.)-(d.)$  The four Gauss diagrams of an oriented loop $h$ with two self-intersection points.}\label{gaussdiagrams.fig}
\end{figure}

Throughout this section, we will use the following facts, which hold for a compact surface $F$ with negative sectional curvature (though compactness is not needed for $(3)$). 
\textit{\begin{enumerate}
\item Nontrivial abelian subgroups of $\pi_1(F)$ are infinite cyclic.
\item There is a unique, maximal infinite cyclic group containing each nontrivial $\alpha \in \pi_1(F)$.
\item Two distinct geodesic arcs with common endpoints cannot be homotopic.
\item Each nontrivial $\alpha\in \fh(F)$ contains a geodesic representative which is unique up to choice of parametrization.
\end{enumerate}}
\noindent The first fact holds by Preissman's Theorem.  The second fact is true if $\partial F\neq\emptyset$ because $\pi_1(F)$ is free.  If $F$ is closed, the second fact follows from the proof of Preissman's Theorem \cite{Chernov}.  The third and fourth facts can be found in \cite{Buser}, as Theorems 1.5.3 and 1.6.6 respectively.\\\\
We now show that for any free homotopy class $\alpha$ on a compact surface, it is possible to choose a representative of $\alpha$ such that no two terms coming from Type 1 intersection points cancel.  This proof is based on ideas in \cite{TuraevViro} and \cite{Chernov}.  Later we will see that if $F=S^2, T^2,$ or the annulus $A$,  geodesic loop on $F$ has no Type 1 self-intersection points, so in Theorem \ref{mainlemma}, we only consider surfaces of negative curvature. 
\begin{thm}\label{mainlemma} Let $F$ be a compact surface equipped with a metric of negative curvature.  Let $\alpha\in \fh(F)$.  If $g$ is a geodesic representative of $\alpha$, then no two terms of $\mu(\alpha)$ corresponding to Type 1 intersection points of $g$ cancel.
\end{thm}
\noindent \textit{Proof.} Throughout this proof, $[\cdot]$ denotes a free homotopy class (either of a geometrical chord diagram or a loop), $[\cdot]_p$ denotes a homotopy class in $\pi_1(F,p)$, and $[\cdot]^p_q$ denotes the homotopy class of a path from $p$ to $q$ with fixed endpoints.  When we concatenate two paths $p_1$ and $p_2$, we write $p_1p_2$, where the path written on the left is the path we traverse first.  \\\\
We write $g=h^n$ for some geodesic loop $h$ and some $n\geq 1$, where $h$ is not a power of another loop.  Suppose $h$ has $m$ self-intersection points, and let $g'$ be a perturbation of $g$ with $mn^2$ Type 1 self-intersection points and $n-1$ Type 2 self-intersection points.  Let $\{p_{i,j}: 1\leq i,j\leq n\}$ and $\{q_{k,l}: 1\leq k,l\leq n\}$ be the sets of $n^2$ self-intersection points corresponding to the (transverse) self-intersection points $p$ and $q$ of $h$ respectively, with the indexing as defined in the previous section.  We assume $[h]$ is nontrivial, since the theorem clearly holds when $[h]$ is trivial ($\mathcal{SI}_0$ is in fact empty). \\\\
We wish to show that the terms of $\mu$ corresponding to points $p_{i,j}$ and $q_{k,l}$ cannot cancel.  We suppose these terms cancel, and derive a contradiction.\\\\
First, we consider the case where $p=q=h(t_1)=h(t_2)$ for $(t_1,t_2)\in \mathcal{SI}_0$, but $i$ and $k$ may or may not be equal, and $j$ and $l$ may or may not be equal.  In other words, $p_{i,j}$ and $q_{k,l}=p_{k,l}$ come from the same set of $n^2$ type 1 self-intersection points.  Let $X=h^1_p$ and let $Y=h^2_p$.  If either $i\neq k$ or $j\neq l$, then by Lemma \ref{type1lemma}, the terms corresponding to $p_{i,j}$ and $p_{k,l}$ are
 $$[(XY)^I X\bullet_p (YX)^JY]-[(YX)^JY\bullet_p (XY)^IX], \text{ and}$$
  $$[(XY)^K X\bullet_p (YX)^LY]-[(YX)^LY\bullet_p (XY)^KX],$$
  for integers $I, J, K, L$ such that $I+J=K+L=n-1$.  If $i=k$ and $k=l$, then $p_{i,j}=p_{k,l}$ corresponds to a single element of $\mathcal{SI}_0$, so we have just the first two of the above terms.  In either case, it suffices to assume that the terms $[(XY)^I X\bullet_p (YX)^JY]$ and $[(YX)^LY\bullet_p (XY)^KX]$ cancel, where $I$ and $K$ may or may not be equal, and $J$ and $L$ may or may not be equal.\\\\
Suppose that
 $$[(XY)^I X\bullet_p (YX)^J Y]=[(YX)^L Y \bullet_p (XY)^K X].$$
 Then there exists $\gamma \in \pi_1(F,p)$ such that 
 \begin{equation}\label{powerssamept1}
 \gamma [(XY)^I X]_p \gamma^{-1}=[(YX)^L Y]_p
 \end{equation}
 and
 \begin{equation}\label{powerssamept2}
 \gamma [(YX)^J Y]_p \gamma^{-1}=[(XY)^K X]_p.
 \end{equation}
 We multiply Equations \ref{powerssamept1} and \ref{powerssamept2} in both possible orders to obtain the equations
 \begin{equation}\label{powerssamept12}
 \gamma [(XY)^I X]_p[(YX)^J Y]_p \gamma^{-1}=[(YX)^L Y]_p[(XY)^K X]_p
 \end{equation}
 and
 \begin{equation}\label{powerssamept21}
 \gamma [(YX)^J Y]_p [(XY)^I X]_p \gamma^{-1}=[(XY)^K X]_p[(YX)^L Y]_p. 
 \end{equation}
 Conjugating Equation \ref{powerssamept12} by $[X]_p\in \pi_1(F,p)$ tells us that $[X]_p\gamma$ and $[(XY)^n]_p$ commute, since $I+J+1=K+L+1=n$.  Similarly, conjugating Equation \ref{powerssamept21} by $[Y]_p$ tells us $[Y]_p\gamma$ and $[(YX)^n]_p$ commute.  Therefore the subgroups $\langle [X]_p\gamma, [(XY)^n]_p \rangle$ and $\langle [Y]_p\gamma, [(YX)^n]_p\rangle$ are infinite cyclic, and are generated by elements $s$ and $t$ of $\pi_1(F,p)$, respectively.  Note that these subgroups are nontrivial since $h$ is nontrivial.  Fact $(2)$ states that each nontrivial element of $\pi_1$ is contained in a unique, maximal infinite cyclic group.  Let $m_1$ and $m_2$ be the generators of the unique maximal infinite cyclic groups containing $[(XY)^n]_p$ and $[(YX)^n]_p$ respectively.  Since $[h]=[XY]=[YX]$ is not freely homotopic to a power of another class, we have that $\langle m_1 \rangle = \langle [XY]_p \rangle$ and $\langle m_2 \rangle=\langle [YX]_p \rangle$.  But $\langle s \rangle$ and $\langle t \rangle$ are also infinite cyclic groups containing $[(XY)^n]_p$ and  $[(YX)^n]_p$, respectively.  By the maximality of the $\langle m_i \rangle$, we have that $\langle s \rangle \leq \langle m_1 \rangle$ and $\langle t \rangle \leq \langle m_2 \rangle$.   This tells us $[X]_p\gamma$ and $[Y]_p\gamma$ are powers of $[XY]_p$ and $[YX]_p$ respectively, so
 \begin{equation}\label{powerssamept2geos}
 \gamma=[X]_p^{-1}([XY]_p)^u=[Y]_p^{-1}([YX]_p)^v.
 \end{equation}
 The powers $u$ and $v$ can be either zero, positive, or negative.  Once we make all possible cancellations in Equation \ref{powerssamept2geos}, we will have two geodesic lassos (one on each side of the equation) formed by products of $X$, $Y$, or their inverses, representing the same homotopy class in $\pi_1(F,p)$.  Therefore these geodesic lassos must coincide.  The geodesic on the left hand side of Equation \ref{powerssamept2geos} can begin by going along either  $X^{-1}$ or $Y$ (depending on the sign of $u$), while the geodesic on the right hand side can begin along either $Y^{-1}$ or $X$ (depending on the sign of $v$).  Therefore $[X]_p$ and $[Y]_p$ must either be powers of the same loop, which is impossible, because we assumed $[h]$ is not a power of another class, or $[X]_p$ and $[Y]_p$ must be trivial, which is impossible because of the definition of $\mathcal{SI}_0$.   Therefore the terms of $\mu$ corresponding to $p_{i,j}$ and $q_{k,l}$ cannot cancel when $p=q=h(t_1)=h(t_2)$.\\\\
Now we will show that the terms of $\mu$ which correspond to $p_{i,j}$ and $q_{k,l}$ cannot cancel when $p$ and $q$ correspond to different ordered pairs in $\mathcal{SI}_0$.  Let $X=h^1_p$, $Y=h^2_p$, $Z=h^1_q$, and $W=h^2_q$.  By Lemma \ref{type1lemma} the terms which $p_{i,j}$ and $q_{k,l}$ contribute to $\mu$ are:
$$[(XY)^I X\bullet_p (YX)^J Y]-[(YX)^J Y \bullet_p (XY)^I X]$$
and
$$[(ZW)^K Z \bullet_q (WZ)^L W]-[(WZ)^L W \bullet_q (ZW)^K Z],$$
where $I+J=K+L=n-1$.  We will suppose that $[(XY)^I X\bullet_p (YX)^J Y]=[(WZ)^L W \bullet_q (ZW)^K Z]$, and derive a contradiction.  Switching the orders of the two loops on both sides of the equation gives us the equality $[(YX)^J Y \bullet_p (XY)^I X]=[(ZW)^K Z \bullet_q (WZ)^L W]$, so if we assume that one of these equalities holds, all four terms above will cancel.  \\\\
As in the case where $p=q$, we will use the equality $[(XY)^I X\bullet_p (YX)^J Y]=[(WZ)^L W \bullet_q (ZW)^K Z]$ to find abelian subgroups of $\pi_1(F,q)$.  To do this, we examine the Gauss diagram of $h$ with two oriented chords corresponding to the self-intersection points $p$ and $q$.  The four possible Gauss diagrams with two oriented chords are pictured in Figure \ref{gaussdiagrams.fig}.   
We use the convention that each oriented chord points from the second branch of $h$ to the first branch of $h$, where the branches of $h$ at a self-intersection point are ordered according to the orientation of $F$.  As shown in Figure \ref{gaussdiagrams.fig}, we let $a_i$, $i=1,..,4$, denote the arcs between the preimages of $p$ and $q$.  We let $b_i$ denote the image of the arc $a_i$ under $h$. \\\\
We first change the basepoint of the first term from $p$ to $q$, replacing $[(XY)^I X\bullet_p (YX)^J Y]$ by $[b_2^{-1}(XY)^I X b_2 \bullet_q b_2^{-1}(YX)^J Y b_2]$. Assuming the terms cancel, we can find $\gamma\in \pi_1(F,q)$ such that
\begin{equation}\label{differentpts1}
\gamma[(WZ)^LW]_q\gamma^{-1}= [b_2^{-1}(XY)^I X b_2]_q
\end{equation}
and
\begin{equation}\label{differentpts2}
\gamma [(ZW)^K Z]_q \gamma^{-1} = [b_2^{-1} (YX)^J Y b_2]_q. 
\end{equation}
Multiplying Equations \ref{differentpts1} and \ref{differentpts2} in both possible orders, and using the fact that $n-1=I+J=K+L$, we have:
\begin{equation}\label{differentpts12}
\gamma[(WZ)^n]_q \gamma^{-1}=[b_2^{-1} (XY)^n b_2]_q
\end{equation}
and
\begin{equation}\label{differentpts21}
\gamma[(ZW)^n]_q \gamma^{-1}=[b_2^{-1} (YX)^n b_2]_q.
\end{equation}
Table \ref{termsforgaussdiagrams} lists the values of $X, Y, Z$ and $W$ in terms of the $b_i$ for each Gauss diagram in Figure \ref{gaussdiagrams.fig}.  
\begin{table}[htdp]
\caption{The values of $X, Y, Z$ and $W$ for each Gauss diagram in Figure \ref{gaussdiagrams.fig}.}
\begin{center}
\begin{tabular}{|c|c|c|c|c|}\hline
{\bf Gauss Diagram} & {\bf $X=h^1_p$} & {\bf $Y=h^2_p$} & {\bf $Z=h^1_q$} & {\bf $W=h^2_q$}\\ \hline
$(a.)$ & $b_2b_3b_4$ & $b_1$ & $b_3$ & $b_4b_1b_2$\\ \hline
$(b.)$ & $b_2b_3b_4$ & $b_1$ & $b_4b_1b_2$ & $b_3$\\ \hline
$(c.)$ & $b_1$ & $b_2b_3b_4$ & $b_3$ & $b_4b_1b_2$\\ \hline
$(d.)$ & $b_2b_3$ & $b_4b_1$ & $b_3b_4$ & $b_1b_2$\\ \hline
\end{tabular}
\end{center}
\label{termsforgaussdiagrams}
\end{table}
This allows us to rewrite Equations \ref{differentpts12} and \ref{differentpts21} just in terms of $\gamma$ and the $b_i$.  Note that for diagrams $(b.)$ and $(c.)$,  we get the same two equations from \ref{differentpts12} and \ref{differentpts21}, because the values of $X$ and $Y$, and the values of $Z$ and $W$, are interchanged.  Therefore it suffices to consider diagrams $(a.)$, $(b.)$, and $(d.)$.  The arguments for diagrams $(a.)$ and $(b.)$ are similar, so we will only examine $(a.)$ and $(d.)$.\\\\
\noindent\textbf{Diagram $(a.)$:}  In this case, $X=b_2b_3b_4$, $Y=b_1$, $Z=b_3$, and $W=b_4b_1b_2$ (see Table \ref{termsforgaussdiagrams}).  We express Equations \ref{differentpts12} and \ref{differentpts21} in terms of the $b_i$ to obtain
\begin{equation}\label{differentpts12bi}
\gamma[(b_4b_1b_2b_3)^n]_q \gamma^{-1}=[b_2^{-1} (b_2b_3b_4b_1)^n b_2]_q
\end{equation}
and
\begin{equation}\label{differentpts21bi}
\gamma[(b_3b_4b_1b_2)^n]_q \gamma^{-1}=[b_2^{-1} (b_1b_2b_3b_4)^n b_2]_q.
\end{equation}
We conjugate Equation \ref{differentpts12bi} by $[b_3]_q^{-1} \in \pi_1(F,q)$ and Equation \ref{differentpts21bi} by $[b_3b_4b_2]_q\in \pi_1(F,q)$ to obtain the equations
$$[b_3]_q^{-1}\gamma[(b_4b_1b_2b_3)^n]_q \gamma^{-1} [b_3]_q = [(b_4b_1b_2b_3)^n]_q$$
and
$$[b_3b_4b_2]_q \gamma [(b_3b_4b_1b_2)^n]_q \gamma^{-1} [b_3b_4b_2]_q^{-1}=[(b_3b_4b_1b_2)^n]_q.$$
Therefore $[b_3]_q^{-1}\gamma$ and $[(b_4b_1b_2b_3)^n]_q$ commute, as do $[b_3b_4b_2]_q \gamma$ and $[(b_3b_4b_1b_2)^n]_q$.  Since abelian subgroups of $\pi_1(F,q)$ are infinite cyclic, the subgroups $\langle [b_3]_q^{-1}\gamma,  [(b_4b_1b_2b_3)^n]_q\rangle$ and $\langle [b_3b_4b_2]_q \gamma, [(b_3b_4b_1b_2)^n]_q\rangle$ are generated by elements $s$ and $t$ in $\pi_1(F,q)$ respectively.  Each nontrivial element of $\pi_1(F,q)$ is contained in a unique, maximal infinite cyclic group by Fact $(2)$.  Let $m_1$ and $m_2$ be the generators of the unique maximal infinite cyclic groups containing $[(b_4b_1b_2b_3)^n]_q$ and $[(b_3b_4b_1b_2)^n]_q$ respectively.  By assumption, $[h]=[b_4b_1b_2b_3]=[b_3b_4b_1b_2]$ is not freely homotopic to a power of another class.  Therefore $\langle m_1 \rangle = \langle [b_4b_1b_2b_3]_q \rangle$ and $\langle m_2 \rangle = \langle [b_3b_4b_1b_2]_q \rangle$.  But $\langle s \rangle$ and $\langle t \rangle$ are also infinite cyclic groups containing $[(b_4b_1b_2b_3)^n]_q$ and $[(b_3b_4b_1b_2)^n]_q$, respectively, so by the maximality of $\langle m_1 \rangle$ and $\langle m_2 \rangle$, we have $\langle s \rangle \leq \langle m_1 \rangle$ and $\langle t \rangle \leq \langle m_2 \rangle$.  Thus $[b_3]_q^{-1}\gamma=([b_4b_1b_2b_3]_q)^u$ and $[b_3b_4b_2]_q\gamma=([b_3b_4b_1b_2]_q)^v$ for some $u$ and $v\in \Z$.  Now 
$$\gamma=[b_3(b_4b_1b_2b_3)^u]_q=[b_2^{-1}b_4^{-1}b_3^{-1}(b_3b_4b_1b_2)^v]_q,$$
so the path homotopy classes $[b_2b_3(b_4b_1b_2b_3)^u]^p_q$ and $[b_4^{-1}b_3^{-1}(b_3b_4b_1b_2)^v]^p_q$ are equal.  Once we cancel $b_i$ with $b_i^{-1}$ wherever possible, $p_1=b_2b_3(b_4b_1b_2b_3)^u$ and $p_2=b_4^{-1}b_3^{-1}(b_3b_4b_1b_2)^v$ will be two geodesic arcs from $p$ to $q$ representing the same path homotopy class.  Therefore $p_1$ and $p_2$ must coincide.  Note that some of the $b_i$ may be trivial.  We know $b_1$ and $b_3$ cannot be trivial because of the definition of $\mathcal{SI}_0$.  Given that $b_2$ or $b_4$ may be trivial, and that $u$ and $v$ may be positive, negative, or zero, we see that $p_1$ can begin along $b_2$, $b_3$, or $b_1^{-1}$ and $p_2$ can begin along $b_4^{-1}$, $b_3^{-1}$, or $b_1$.  Thus $p_1$ and $p_2$ can only coincide if the beginnings of the arcs $b_i$ and $b_j^{\pm 1}$, as well as the initial velocity vectors of these arcs, coincide for some $i\neq j$.  This is impossible since $h$ is a geodesic which is not homotopic to a power of another loop. \\\\

\noindent\textbf{Diagram $(d.)$:} 
 In this case, $X=b_2b_3$, $Y=b_4b_1$, $Z=b_3b_4$, and $W=b_1b_2$; see Table \ref{termsforgaussdiagrams}.  We rewrite Equations \ref{differentpts12} and \ref{differentpts21} in terms of the $b_i$ to obtain
 \begin{equation}\label{differentpts12bid}
\gamma[(b_1b_2b_3b_4)^n]_q \gamma^{-1}=[b_2^{-1} (b_2b_3b_4b_1)^n b_2]_q
\end{equation}
and
\begin{equation}\label{differentpts21bid}
\gamma[(b_3b_4b_1b_2)^n]_q \gamma^{-1}=[b_2^{-1} (b_4b_1b_2b_3)^n b_2]_q.
\end{equation}
We conjugate Equation \ref{differentpts12bid} by $[b_1b_2]_q \in \pi_1(F,q)$ and Equation \ref{differentpts21bid} by $[b_3b_2]_q\in \pi_1(F,q)$ to obtain the equations
$$[b_1b_2]_q\gamma[(b_1b_2b_3b_4)^n]_q \gamma^{-1} [b_1b_2]_q^{-1} = [(b_1b_2b_3b_4)^n]_q$$
and
$$[b_3b_2]_q \gamma [(b_3b_4b_1b_2)^n]_q \gamma^{-1} [b_3b_2]_q^{-1}=[(b_3b_4b_1b_2)^n]_q.$$
Therefore $[b_1b_2]_q\gamma$ and $[(b_1b_2b_3b_4)^n]_q$ commute, as do $[b_3b_2]_q \gamma$ and $[(b_3b_4b_1b_2)^n]_q$.  Since abelian subgroups of $\pi_1(F,q)$ are infinite cyclic, the subgroups $\langle [b_1b_2]_q\gamma,  [(b_1b_2b_3b_4)^n]_q\rangle$ and $\langle [b_3b_2]_q \gamma, [(b_3b_4b_1b_2)^n]_q\rangle$ are generated by elements $s$ and $t$ in $\pi_1(F,q)$ respectively.  Each nontrivial element of $\pi_1(F,q)$ is contained in a unique, maximal infinite cyclic group by Fact $(2)$.  Let $m_1$  and $m_2$ be the generators of the unique, maximal infinite cyclic groups containing $[(b_1b_2b_3b_4)^n]_q$ and $[(b_3b_4b_1b_2)^n]_q$, respectively.  Since $[h]=[b_1b_2b_3b_4]=[b_3b_4b_1b_2]$ is not freely homotopic to a power of another class, we have $\langle m_1 \rangle = \langle [b_1b_2b_3b_4]_q \rangle$ and $\langle m_2 \rangle = \langle[b_3b_4b_1b_2]_q \rangle$.  But $\langle s \rangle$ and $\langle t \rangle$ are also infinite cyclic groups containing $[(b_1b_2b_3b_4)^n]_q$ and $[(b_3b_4b_1b_2)^n]_q$, respectively.  Thus by the maximality of the $\langle m_i \rangle$, we have $\langle s \rangle \leq \langle m_1 \rangle$ and $\langle t \rangle \leq \langle m_2 \rangle$.  Hence $[b_1b_2]_q\gamma=([b_1b_2b_3b_4]_q)^u$ and $[b_3b_2]_q\gamma=([b_3b_4b_1b_2]_q)^v$ for some $u$ and $v\in \Z$.  Now 
$$\gamma=[b_2^{-1}b_1^{-1}(b_1b_2b_3b_4)^u]_q=[b_2^{-1}b_3^{-1}(b_3b_4b_1b_2)^v]_q,$$
so the path homotopy classes $[b_1^{-1}(b_1b_2b_3b_4)^u]^p_q$ and $[b_3^{-1}(b_3b_4b_1b_2)^v]^p_q$ are equal.  Once we cancel $b_i$ with $b_i^{-1}$ wherever possible, $p_1=b_1^{-1}(b_1b_2b_3b_4)^u$ and $p_2=b_3^{-1}(b_3b_4b_1b_2)^v$ will be two geodesic arcs from $p$ to $q$ representing the same path homotopy class. Therefore $p_1$ and $p_2$ must coincide.  Again, some of the $b_i$ may be trivial.  Because of the definition of $\mathcal{SI}_0$, adjacent arcs (e.g. $b_2$ and $b_3$ or $b_4$ and $b_1$) cannot both be trivial.  If arcs opposite each other (e.g. $b_2$ and $b_4$) are both trivial, then $p=q$; that case was already examined.  So we may assume at most one of the $b_i$ is trivial.  Depending on whether $u$ is positive, negative, or zero, and on which $b_i$ is trivial, $p_1$ can begin along either $b_1^{-1}$, $b_2$, $b_4^{-1}$, $b_3$, or $p_1$ can be trivial (if $u=0$ and $b_1$ is trivial).  Similarly, $p_2$ can begin along $b_4$, $b_1$, $b_3^{-1}$, or $b_2^{-1}$, or $p_2$ can be trivial (if $v=0$ and $b_3$ is trivial).  Therefore, in order for the $p_i$ to coincide, either the beginnings and initial velocity vectors of the arcs $b_i$ and $b_j^{\pm 1}$ must coincide for some $i\neq j$, which is impossible since $h$ is a geodesic and is not a power of another loop, or both $p_i$ must be trivial.  But if both $p_i$ are trivial, then $b_1$ and $b_3$ are both trivial, and we assumed at most one of the $b_i$ are trivial, so this is impossible as well.\qed \\\\
The following lemma allows us to reduce the proofs of Theorems \ref{main2} and \ref{nonprimthm} to the case where $F$ is compact.
\begin{lem} \label{compact} Suppose $F$ is noncompact, and let $g:S^1\rightarrow F$.  Suppose $t(\mu([g]))=T$, where $[g]\in\fh(F)$.  Then there exists a compact subsurface $F_C$ of $F$ containing $\im(g)$ such that $t(\mu([g]_C))=T$, where $[g]_C$ is the class of $g$ in $\fh(F_C)$. \end{lem}
\noindent\textit{Proof.}  First note that we may assume $g$ has finitely many transverse self-intersection points (if it did not, we could perturb $g$ slightly to obtain a loop that does).  Suppose that the pair of terms $[g^1_p\bullet_p g^2_p]$ and $[g^2_q\bullet_q g^2_q]$ of $\mu([g])$ cancel.  Let $D$ be a chord diagram with one oriented chord attached to one copy of $S^1$ at its endpoints.  Let $D_p: D\rightarrow F$ and $D_q:D\rightarrow F$ be the geometrical chord diagrams associated with $[g^1_p \bullet_p g^2_p]$ and $[g^2_q \bullet_q g^1_q]$ respectively.  Since  $[g^1_p \bullet_p g^2_p]$ and $[g^2_q \bullet_q g^1_q]$ cancel, we have a homotopy $H_{p,q}:D\times [0,1] \rightarrow F$ between $D_p$ and $D_q$.  Since $\im(H_{p,q})$ is compact, and since $g$ has finitely many self-intersection points, we may choose a compact subsurface $F_C$ of $F$ containing $\im(H_{p,q})$ for all pairs $(p,q)$ corresponding to terms that cancel.  Note that once $F_C$ contains $\im(H_{p,q})$, $F_C$ must also contain the image of $g$.  If we compute $\mu([g]_C)$ on $F_C$, the terms  $[g^1_p \bullet_p g^2_p]$ and $[g^2_q \bullet_q g^1_q]$ will cancel, since $H_{p,q}$ can be viewed as a homotopy in $F_C$.  Thus $t(\mu([g]))\geq t(\mu([g]_C))$.  Furthermore, any terms which cancel on $F_C$ must cancel on $F$, so the inequality becomes an equality.
\qed\\\\
Now we state our main results.   Theorems \ref{main2} and \ref{nonprimthm} are stated as Corollaries \ref{cor1} and \ref{cor2} of Theorem \ref{mainlemma}, respectively, though for now we restrict Theorem \ref{nonprimthm} to the case where $\alpha$ is not a power of another class.  Recall that Theorem \ref{main2} states that $\mu(\alpha)$ is zero if and only if $\alpha$ is a power of a simple class, and Theorem \ref{nonprimthm} gives a formula for the minimal self-intersection number of $\alpha$.

\begin{cor}\label{cor1} Let $\alpha \in \fh(F)$.  Then $\mu(\alpha)=0$ if and only if $\alpha=\beta^n$ for some $\beta\in\pi_1(F)$, where $m(\beta)=0$.
\end{cor}
\noindent\textit{Proof.}  Suppose $\alpha=\beta^n$ for some $\beta\in\pi_1(F)$, where $m(\beta)=0$.  We may assume $n>0$.  We will compute $\mu(\alpha)$, and see that $\mu(\alpha)=0$.  We begin by choosing a simple representative $h$ of $\beta$ and a point $p$ on $h$. Then $g=(h_p)^n$ is a representative of $\alpha$.  We perturb $g$ slightly so that it has $n-1$ self-intersection points (of type 2), all with image $p$. Now
$$\mu(\alpha)=\sum_{i=1}^{n-1} [(h_p)^i\bullet_p (h_p)^{n-i}]-[(h_p)^{n-i}\bullet_p (h_p)^i],$$\\
which equals zero, since the positive term corresponding to $i=k$ cancels with the negative term corresponding to $i=n-k$.  Note that this argument actually shows that the terms of $\mu$ corresponding to Type 2 self-intersection points always cancel with each other, even if $m(\beta)\neq 0$.\\\\
To prove the converse, we assume $\alpha$ cannot be written as a power of a simple class.  By Lemma \ref{compact}, we may assume $F$ is compact.  Therefore, either $F=S^2, A,$ or $T^2$, or $F$ can be equipped with a metric of negative curvature.  In this situation, $F$ clearly cannot be $S^2$ or $A$.\\\\
$F$ also cannot be $T^2$.  Let $g$ be a representative of $\alpha$, and consider $[g]_p\in\pi_1(T^2,p)$ for some $p$ on the image of $g$.  Let $a$ and $b$ be the generators of $\pi_1(T^2,p)$ which lift to the paths from $(0,0)$ to $(1,0)$ and $(0,1)$ respectively under the usual covering map $\mathbb{R}^2\rightarrow \mathbb{R}^2/\Z^2$.  Then we can write $[g]_p=([a^{e_1}b^{e_2}]_p)^k$ where the $e_i$ are relatively prime. Taking $\beta=[a^{e_1}b^{e_2}]_p$, we have $m(\beta)=0$, since the lift of $a^{e_1}b^{e_2}$ is homotopic to a path from $(0,0)$ to $(e_1,e_2)$.\\\\
Now we may assume $F\neq S^2, A$ or $T^2$, so we can apply Theorem \ref{mainlemma}.  We choose a geodesic representative $g$ of $\alpha$, and write $g=h^n$ where $h$ is not a power of another loop in $\pi_1(F)$.  Suppose $h$ has $m$ self-intersection points.  Then we can perturb $g$ to obtain a loop $g'$ with $mn^2+n-1$ self-intersection points, where $mn^2$ of the self-intersections are of Type 1, and $n-1$ are of Type 2.  By Theorem \ref{mainlemma}, no two terms of $\mu$ corresponding to Type 1 self-intersection points can cancel.  As we saw above, all of the Type 2 terms cancel with each other.  So, after all cancellations are made, $2mn^2$ terms of $\mu$ remain.  Since $m\geq 1$, we know $\mu(\alpha)\neq 0$.
\qed
\begin{cor}\label{cor2} Let $\alpha\in\fh(F)$ be primitive.  Then $m(\alpha)=t(\mu(\alpha))/2$.  Thus $\mu$ computes the minimal self intersection number of $\alpha$.
\end{cor}
\noindent\textit{Proof.} By Lemma \ref{compact}, we may assume $F$ is compact.  If $\alpha\neq \beta^n$ for any $\beta \in \pi_1(F)$ and $|n|>1$, we can choose a geodesic representative $g$ of $\alpha$ such that the self-intersection points of its perturbation $g'$ are all of Type 1.  If $F\neq S^2, A$ or $T^2$, then by Theorem \ref{mainlemma}, no two terms of $\mu(\alpha)$ can cancel.  Thus $m(\alpha)=t(\mu(\alpha))/2$. \\\\
 If $F=S^2, T^2$ or $A$,  we can choose a representative of $\alpha$ with no self-intersection points, implying that the equality $m(\alpha)=t(\mu(\alpha))/2$ holds trivially.  This is clear for $S^2$ and $A$.  If $F=T^2$, we choose a representative $g$ of $\alpha$, and consider $[g]_p\in \pi_1(T^2,p)$ for some point $p$ on the image of $g$.  If $a$ and $b$ are generators the of $\pi_1(T^2,p)$ as in the proof of Corollary \ref{cor1}, we can write $[g]_p=([a^{e_1}b^{e_2}]_p)^k$ where the $e_i$ are relatively prime.  Since $\alpha$ is not a power of another class, $k=\pm 1$, so $[g]_p= ([a^{e_1}b^{e_2}]_p)^{\pm 1}$.  This class has a simple representative, as shown in the proof of Corollary \ref{cor1}.   \qed

\subsection{Using $\mu$ to compute the minimal self-intersection number of a class which is not primitive.}  If we combine our results above with those of Hass and Scott in \cite{HassScott}, we can use $\mu$ to compute the minimal self-intersection number of any class $\alpha$, even if $\alpha$ is not primitive.  We will use the following result, which we state for orientable surfaces (their result is more general):
\begin{lem}[Hass, Scott]  Let $f$ be a loop on $F\neq S^2$ in general position.  Suppose that $f$ is a representative of the class $\alpha=\beta^n$, where $\beta$ is a non-trivial and primitive element of $\pi_1(F)$.  Let $\tilde{F}$ be the universal cover of $F$ and let $F_{\beta}$ denote the quotient of $\tilde{F}$ by the cyclic subgroup of $\pi_1(F)$ generated by $\beta$.  Let $f_{\beta}:S^1\rightarrow F_{\beta}$ be the lift of $F$, and let $l$ denote one of the lines in $\tilde{F}$ above $f_{\beta}(S^1)$.  Then $f$ has least possible self-intersection if and only if $f_{\beta}$ has least possible self-intersection and, for all $\gamma$ in $\pi_1(F)$ such that $\gamma$ is not a power of $\beta$, the intersection $\gamma l\cap l$ consists of at most one point. \label{HassScottLemma}
\end{lem}
\noindent We now use Lemma \ref{HassScottLemma} to get the following corollary, which gives a formula for $m(\alpha)$ in terms of $m(\beta)$, where $\alpha$ is a power of $\beta$ and $\beta$ is primitive.  For a similar statement for compact surfaces with boundary, see the remark following Theorem 2 in \cite{Tan}.  (Note that as it seems, the statement in \cite{Tan} contains a typo, namely $ps^2+(p-1)$ should be $ps^2+(s-1)$).  
\begin{cor} Let $\alpha \in \fh(F)$ be a nontrivial class such that $\alpha=\beta^n$ and $\beta$ is a primitive class, where $F\neq S^2$ and $n\geq 1$.  Then $m(\alpha)=n^2m(\beta)+n-1$. \label{corminintpower}
\end{cor}
\noindent\textit{Proof.}  First we remark that one can deduce this fact from the proof of Lemma \ref{HassScottLemma} in \cite{HassScott}, rather than from the statement of the lemma itself.  We choose to do the latter to avoid repeating most of the proof of Lemma \ref{HassScottLemma}. \\\\
Equip $F$ with a flat or hyperbolic metric.  Then $\tilde{F}$ is diffeomorphic to $\mathbb{R}^2$, so $F_\beta$ is a cylinder.  We begin by choosing a geodesic representative $g$ of $\alpha$ and perturb $g$ to obtain a new loop $g'$.  If $n=1$, we take $g'=g$.  Otherwise, we lift $g$ to obtain a geodesic loop $g_{\beta}$ in $F_{\beta}$.  We then take a $C^\infty$-small perturbation of $g_{\beta}$ to obtain a loop $g'_\beta$ with $n-1$ self-intersection points, and we let $g'$ be the projection of $g'_\beta$ back to $F$.  We note that $g'_\beta$ has minimal self-intersection (see, for example, \cite{TuraevViro}).  As mentioned in the discussion leading up to Theorem \ref{mainlemma}, the number of self-intersection points of $g'$ is $n^2m(\beta)+n-1$.  Let $l'$ denote one of the lines above $g'$ in $\tilde{F}$.  To show that $g'$ has minimal self-intersection, it remains to show that for all $\gamma\in \pi_1(F)$ such that $\gamma$ is not a power of $\beta$, the intersection $\gamma l'\cap l'$ consists of at most one point.  Suppose that for some $\gamma$, the lines $\gamma l'$ and $l'$ intersect in two points $p$ and $q$.  These lines bound a 2-gon in $\tilde{F}$.  Let $A_1$ and $A_2$ denote the two arcs of the 2-gon in $\tilde{F}$.  Let $\rho$ denote the covering map from $\tilde{F}$ to $F$.  Since $A_1$ and $A_2$ are homotopic with common ends, their projections $\rho(A_1)$ and $\rho(A_2)$ must be as well.  Furthermore the arcs $\rho(A_1)$ and $\rho(A_2)$ are contained in the image of $g'$.  Since $g'$ is a perturbation of a geodesic, $\rho(p)$ and $\rho(q)$ must be Type 2 self-intersection points of $g'$.  Thus each $\rho(A_i)$ is an arc in the image of $g'$ connecting two Type 2 self-intersection points, and hence must be homotopic to $g^k$ for some $k\neq 0$.  This implies $\gamma$ is a power of $\beta$.  Hence if $\gamma$ is not a power of $\beta$, the intersection of $\gamma l'$ and $l'$ contains at most one point.
\qed\\\\
Recall that, if we make our usual choice of a perturbation of a geodesic $g'\in \alpha$, then two terms of $\mu(\alpha)$ cancel if and only if they correspond to Type 2 self-intersection points.  (The observation that terms corresponding to Type 2 self-intersection points must cancel with other such terms is contained in the proof of Corollary \ref{cor1}.)  Together with Corollary \ref{corminintpower}, this allows us to use $\mu$ to compute the minimal self-intersection number of a class which may not be primitive.
\begin{cor}  Let $\alpha\in \fh(F)$ be a nontrivial class such that $\alpha=\beta^n$ where $\beta$ is primitive and $n\geq 1$.  Then $m(\alpha)=t(\mu(\alpha))/2+n-1$.
\end{cor}
\noindent\textit{Proof.}  By Theorem \ref{mainlemma}, we know $t(\mu(\alpha))/2=m(\beta)n^2$, i.e., the number of Type 1 self-intersection points.  Since $m(\alpha)=m(\beta)n^2+n-1$, we have $m(\alpha)=t(\mu(\alpha))/2+n-1$.
\qed

\section{An Example}  
\begin{figure}[h]\center
\scalebox{1.0}{\includegraphics[width=12cm]{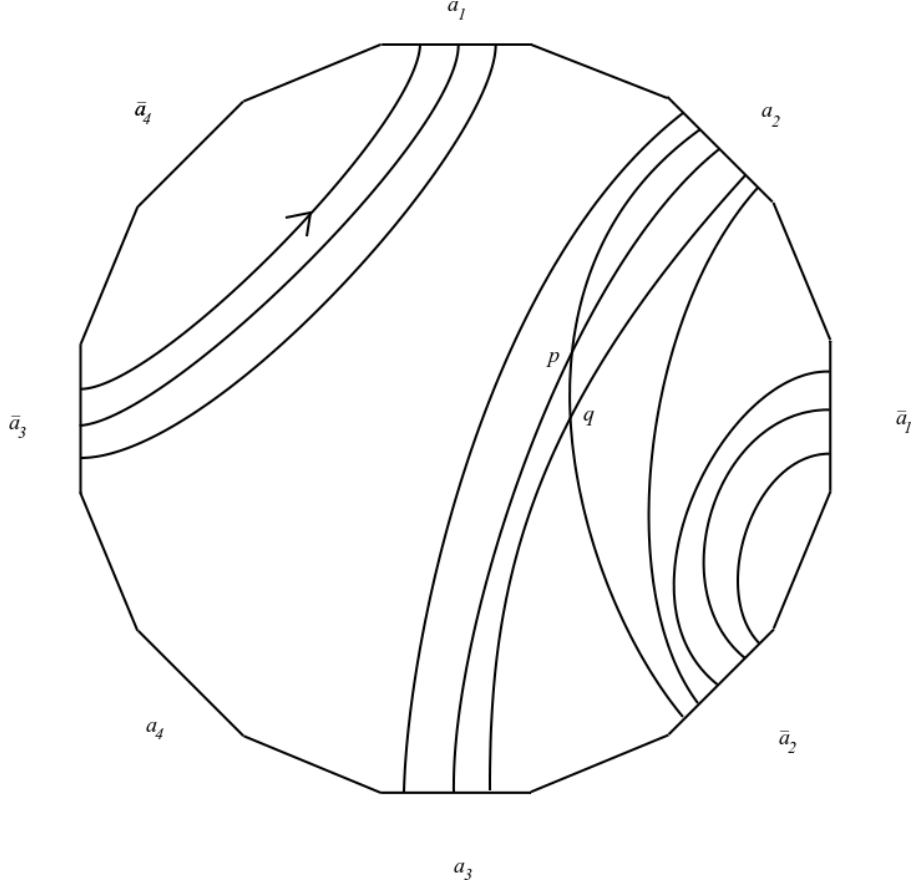}}  
\caption{The class $a_3a_1\bar{a}_2a_3a_1\bar{a}_2a_3a_1\bar{a}_2\bar{a}_2\bar{a}_2$ on the punctured surface of genus two.}\label{genus2puncturedexample.fig}
\end{figure}

The following example illustrates how one can compute the minimal self-intersection number of a free homotopy class algorithmically using $\mu(\alpha)$.  We will use $\mu$ to compute the minimal intersection number of a class 
$$\alpha=a_3a_1\bar{a}_2a_3a_1\bar{a}_2a_3a_1\bar{a}_2\bar{a}_2\bar{a}_2$$
 on the punctured surface of genus two, with surface word $a_1a_2\bar{a}_1\bar{a}_2a_3a_4\bar{a}_3\bar{a_4}.$  This homotopy class has Turaev cobracket zero, as shown in Example 5.8 of \cite{Chas}.    A representative of $\alpha$ with two self-intersection points is pictured in Figure \ref{genus2puncturedexample.fig}.  Using $\mu$, along with the methods in \cite{Chernov}, we show that the minimal self-intersection number of $\alpha$ is 2.
 
We compute $\mu$, and find that
$$\mu(\alpha)=a_3a_1\bar{a}_2a_3a_1\bar{a}_2\bar{a}_2\bullet_p\bar{a}_2a_3a_1\bar{a}_2-\bar{a}_2a_3a_1\bar{a}_2\bullet_p a_3a_1\bar{a}_2a_3a_1\bar{a}_2\bar{a}_2 $$
$$+ a_3a_1\bar{a}_2\bar{a}_2\bullet_q \bar{a}_2a_3a_1\bar{a}_2a_3a_1\bar{a}_2-\bar{a}_2a_3a_1\bar{a}_2a_3a_1\bar{a}_2\bullet_q a_3a_1\bar{a}_2\bar{a}_2.$$
It is easy to see that terms one and two cannot cancel, since there is no $t\in \pi_1(F)$ such that $ta_3a_1\bar{a}_2a_3a_1\bar{a}_2\bar{a}_2 t^{-1}=\bar{a}_2a_3a_1\bar{a}_2$.  In general, one can check whether two reduced words in the generators of $\pi_1$ are in the same conjugacy class by comparing two cyclic lists, since the conjugacy class of a reduced word in a free group consists of all the cyclic permutations of that word. Similarly terms three and four cannot cancel.

Next, we show terms one and four cannot cancel.  First we conjugate term one by $\bar{a}_2$ so that its first loop matches the first loop of term four.  The first term becomes
$$\bar{a}_2a_3a_1\bar{a}_2a_3a_1\bar{a}_2 \bullet_p \bar{a}_2 \bar{a}_2a_3a_1.$$
If terms one and four cancel, we can find $s\in \pi_1(F)$ such that $s$ and $\bar{a}_2a_3a_1\bar{a}_2a_3a_1\bar{a}_2$, and such that 
$$s\bar{a}_2 \bar{a}_2a_3a_1 s^{-1}=a_3a_1\bar{a}_2\bar{a}_2.$$
Abelian subgroups of $\pi_1(F)$ are infinite cyclic, so the subgroup $\langle s, \bar{a}_2a_3a_1\bar{a}_2a_3a_1\bar{a}_2 \rangle$ is generated by some $u\in \pi_1(F)$.  We will see that $\bar{a}_2a_3a_1\bar{a}_2a_3a_1\bar{a}_2=u^{\pm1}$.  Suppose $\bar{a}_2a_3a_1\bar{a}_2a_3a_1\bar{a}_2=u^i$, and consider this relation in the abelianization of $\pi_1(F)$, where the generator $a_j$ of $\pi_1(F)$ is sent to the generator of $\Z^4$ with a 1 in position $j$ and zeroes elsewhere.  Hence this relation becomes $(2,-3,2,0)=i(u_1,u_2,u_3,u_4)$, which only has solutions when $i=\pm 1$.  Therefore $s=(\bar{a}_2a_3a_1\bar{a}_2a_3a_1\bar{a}_2)^k$ for some $k\in\Z$.  However, the relation
$$(\bar{a}_2a_3a_1\bar{a}_2a_3a_1\bar{a}_2)^k\bar{a}_2 \bar{a}_2a_3a_1(\bar{a}_2a_3a_1\bar{a}_2a_3a_1\bar{a}_2)^{-k}=a_3a_1\bar{a}_2\bar{a}_2$$
cannot hold for any value of $k$.  Hence none of the terms of $\mu$ cancel, and $m(\alpha)=2$.\\\\
In general, when deciding whether two terms cancel, we first verify that the first loops in each term are indeed in the same conjugacy class.  If this is the case, we will have two elements $s$ and $t$ that commute ($t=\bar{a}_2a_3a_1\bar{a}_2a_3a_1\bar{a}_2$ in the example above), and we know $\langle s, t\rangle=\langle u \rangle$.  Since every nontrivial element of $\pi_1$ is contained in a unique maximal infinite cyclic subgroup, we can write $\langle s, t \rangle=\langle u\rangle \leq \langle m \rangle$ where $m$ is the generator of the maximal infinite cyclic subgroup containing $u$.  We need to find $m$ given $t$, as our goal is to write $s$ as a power of $m$.  To find $m$ given $t$, we first cyclically reduce $t$ and write the result as a word $t_r$ in the generators of $\pi_1$.  Since $t_r$ is cyclically reduced, then in order to be a power of another element, it must look like a concatenation of $i$ copies of some cyclically reduced word $w$ in the generators of $\pi_1$, so we can determine $i$ and $w$ where $i$ is as large (in absolute value) as possible.  Now $t_r=w^i$, and there exists $c\in \pi_1$ such that $t=ct_rc^{-1}=(cwc^{-1})^i$, so $m=cwc^{-1}$.  We can now write $s$ as a power of $m$ and finish the algorithm as in the example above.

\section{Algebraic Properties of $\mu$}
\begin{figure}[h]\center
\scalebox{1.0}{\includegraphics[width=4cm]{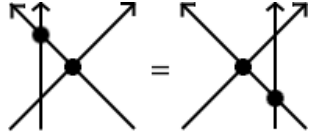}}  
\caption{We can slide one exterior chord past another.}
\label{slidechord.fig}
\end{figure}

\begin{figure}[h]\center
\scalebox{1.0}{\includegraphics[width=4cm]{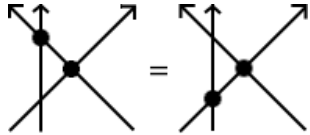}}  
\caption{We can slide one exterior from one core circle to another.}
\label{slidechord2.fig}
\end{figure}

\begin{figure}[h]\center
\scalebox{1.0}{\includegraphics[width=10cm]{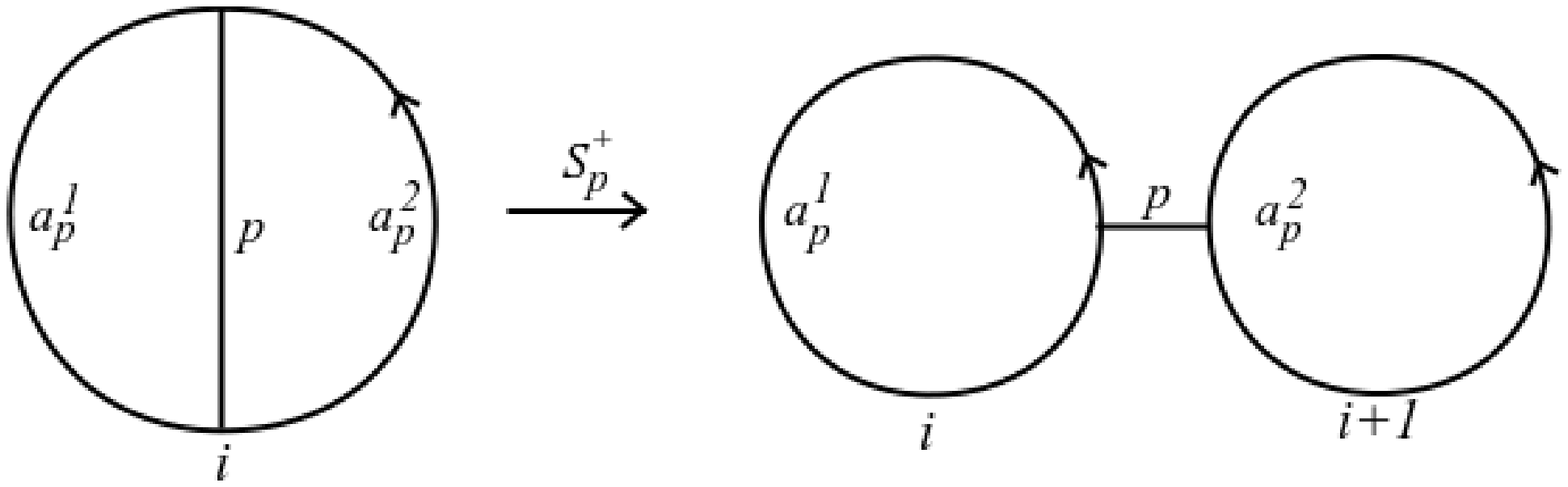}}  
\caption{The splitting map $S_p^+$.}\label{split12.fig}
\end{figure}

\begin{figure}[h]\center
\scalebox{1.0}{\includegraphics[width=10cm]{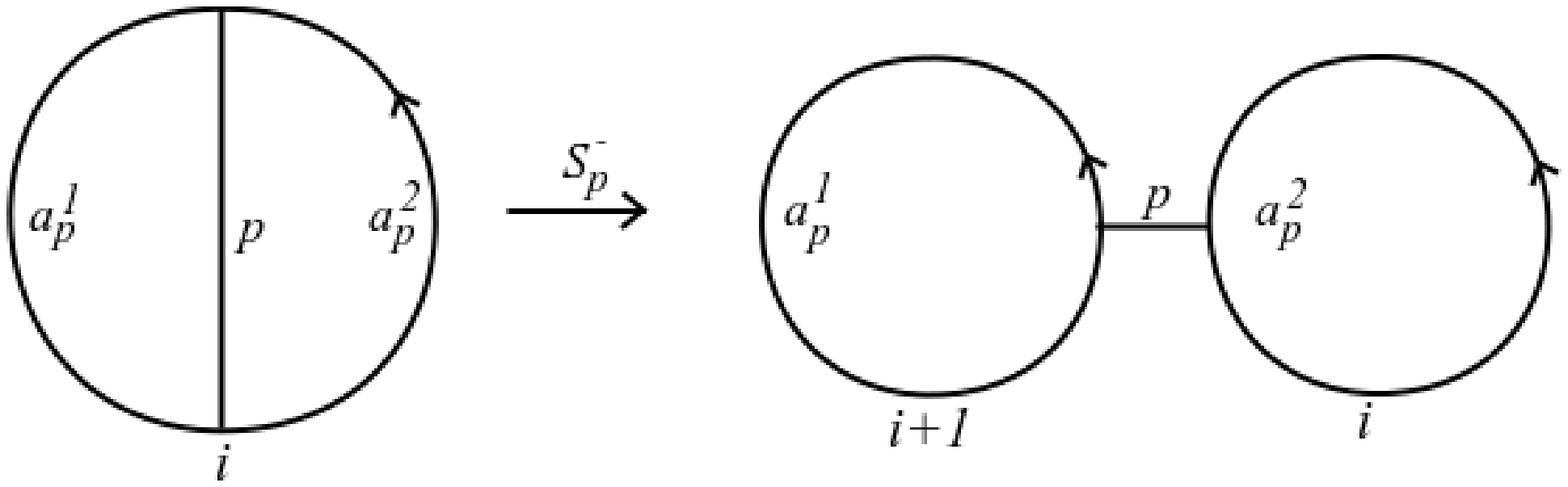}}  
\caption{The splitting map $S_p^-$.}\label{split21.fig}
\end{figure}

We conclude by investigating properties of $\mu$ which allow one to view $\mu$ as a generalization of a Lie cobracket.  In particular, we exhibit analogues of the following properties of the Turaev cobracket $\Delta$:
\begin{itemize}
\item $\Delta$ satisfies co-skew symmetry:  $\tau \circ \Delta = -\Delta$, where $\tau(a \otimes b)=b\otimes a$.
\item $\Delta$ satisfies the co-Jacobi identity: $(1+\omega+\omega^2)\circ Id\otimes \Delta \circ \Delta=0$, and $\omega(a\otimes b\otimes c)= b\otimes c\otimes a$.
\end{itemize}
We begin by modifying the definition of $\mu$ given for free loops, and then extend this operation to certain chord diagrams in the Andersen-Mattes-Reshetikhin algebra.  For the purposes of computing the minimal self-intersection number, this definition is equivalent to the previous one.  However, it is easier to state the analogues of co-skew symmetry and the co-Jacobi identity for this modified definition.\\\\
Recall that $M$ denotes the free $\mathbb{Z}$-module generated by the set of chord diagrams on $F$, and $N$ denotes the submodule generated by the $4T$-relations in Figure \ref{4term.fig} (and relations obtained from them by reversing orientations on arrows).  The Andersen-Mattes-Reshetikhin bracket is defined on the quotient $ch=M/N$.  Given a chord $p$ of a geometrical chord diagram, we say $p$ is an \textit{external chord} if the endpoints of $p$ lie on distinct core circles.  Otherwise, we call $p$ an \textit{internal chord}.  Now we let $M_e$ denote the free $\Z$-module generated by diagrams with only external chords, and let $N_e$ denote the submodule of $M_e$ generated by the relations in Figures \ref{slidechord.fig} and \ref{slidechord2.fig}, their mirror images, and relations obtained from these by reversing the orientation on any branch in the picture.  From these relations, one can obtain the $4T$ relations of the Andersen-Mattes-Reshetikhin algebra in which the orientations of the arcs in the four pictures are identical. We will define $\mu$ on $ch_e=M_e/N_e$.  From now on, we also assume that the core circles of our chord diagrams are labeled with the digits $\{1,...,n\}$, where $n$ is the number of core circles in the diagram.\\\\
The operation $\mu$ will add a chord between the preimages of each self-intersection point of each core circle, in the same way as before.  However, rather than adding an oriented chord whose orientation induces a labeling on the resulting two loops in the image of the geometrical chord diagram, we split the abstract diagram into two labeled loops connected by an unoriented chord.\\\\
We define two maps $S_p^+$ and $S_p^-$ which split a core circle labelled $i$ into two new core circles labeled $i$ and $i+1$ (see Figures \ref{split12.fig} and \ref{split21.fig}).  Suppose $a:C_i\rightarrow F$ is the restriction of our original chord diagram to the $i^{th}$ core circle.  The map $S_p^+$ (respectively $S_p^-$) maps the new circle labeled $i$ (respectively $i+1$) to $a^1_p$ and the new circle labeled $i+1$ (respectively, $i$) to $a^2_p$.  The splitting map also increases by one the label on all core circles formerly labeled with numbers greater than or equal to $i+1$.  Note that the image of smooth chord diagram under $S_p^\pm$ may not be smooth at $p$, but we can always find a smooth diagram in its homotopy class.  \\\\

\subsection{The Definitions of $\mu_i$ and $\mu$.}
Now let $\mathcal{SI}^i_0$ denote the set of self-intersection points of the $i^{th}$ core circle in our chord diagram such that neither of the maps formed by the splitting map at that point are homotopy trivial.  Let $D_p$ denote the diagram obtained by adding a chord between $t_1$ and $t_2$ in a chord diagram $D$ where $(t_1,t_2)\in \mathcal{SI}^i_0$ and $p=D(t_1)=D(t_2)$.  We define 
$$\mu_i(\mathcal{D})=\sum_{p\in\mathcal{SI}^i_0} [S^+_p(D_p)]-[S^-_p(D_p)],$$
and we let 
$$\mu(\mathcal{D})=\sum_{i=1}^n \mu_i(\mathcal{D}).$$
Using linearity, we extend this definition so that $\mu$ becomes a map from $ch_e$ to $ch_e$.\\\\
\textit{Remark:} The identification in Figure \ref{slidechord2.fig} arises because $\mu$ splits the original core circle into two new circles.  If we define $\mu$ without the splitting map, this identification is not needed.
\subsection{The maps $\mu$ and $\mu_i$ are well-defined.}
We need to verify that $\mu_i(\mathcal{D})$ does not depend on the choice of diagram $D$ in $\mathcal{D}$.  Clearly $\mu_i$ does not change when $D$ undergoes regular isotopy.  Applying the first Reidemeister move to $D$ does not affect $\mu_i([D])$ because we only sum over self-intersection points such that the two new maps produced by $S^\pm_p$ are not homotopy trivial.  The fact that $\mu_i$ does not change under other elementary moves follows from either other elementary moves or identifications we make:
\begin{itemize}
\item {\bf Second Reidemeister Move:} This follows from the move in Figure \ref{type2movewithdot.fig}.
\item {\bf Third Reidemeister Move:} This follows from the move in Figure \ref{type3movewithdot.fig}.
\item {\bf The move in Figure \ref{type2movewithdot.fig}:}  Because the existing chord is an exterior chord, we do not sum over the intersection point in this diagram because it must be an intersection point of different core circles. 
\item {\bf The move in Figure  \ref{type3movewithdot.fig}:}  This follows from the following identifications in Figures \ref{slidechord.fig} and \ref{slidechord2.fig}.  Note that the $4T$ relations in the Andersen-Mattes-Reshetikhin algebra imply that their bracket is invariant under this move.  In our case, the relations only contain two terms because the fact that the chord is an exterior chord implies we only sum over at most one of the intersection points before and after the move in Figure \ref{type3movewithdot.fig}.
\item {\bf The identification in Figure \ref{slidechord.fig}:}  This follows from the following identifications in Figures \ref{slidechord.fig} and \ref{slidechord2.fig}:
\item {\bf The identification in Figure \ref{slidechord2.fig}:} Because the chords in the picture are exterior chords, the one intersection point cannot be a self-intersection point, so we do not sum over it.
\end{itemize}
\subsection{Algebraic Properties of $\mu$ and $\mu_i$.}
Let $\tau_i$ be a map which swaps labels the $i$ and $i+1$ in a chord diagram with enumerated core circles.  We have the following analogue of co-skew symmetry:
\begin{prop} $\tau_i \circ \mu_i=-\mu_i.$
\end{prop}
\noindent\textit{Proof.} Clear.
\qed\\\\
Let $\omega_i$ be a map which cyclically permutes the labels $i,i+1, i+2$.  Specifically, $\omega_i$ decreases the labels $i+1$ and $i+2$ by one, and sends $i$ to $i+2$.  We have the following analogue of the co-Jacobi identity:
\begin{prop} For all $1\leq i \leq n-1$, $(1+\omega_i +\omega_i^2)\circ \mu_{i+1}\circ \mu_i=0.$ \label{generalcojacobi}
\end{prop}
\noindent\textit{Proof.} It suffices to show that, for a diagram with one core circle $C_1$, we have $(1+\omega_1 +\omega_1^2)\circ \mu_{2}\circ \mu_1=0.$  Each diagram in the sum $\mu_2\circ\mu_1(C_1)$ has three core circles and two chords corresponding to two self-intersection points of $C_1$, which we call $p$ and $q$.  The idea of the proof is as follows:  $\mu_2\circ\mu_1(C_1)$ is a sum of four terms, two of which are positive and two of which are negative.  Relative to some fixed initial order on the three core circles, the labelings on each diagram form an even or odd permutation in $S_3$, the symmetric group on three elements.  The labelings on two of the four terms of $\mu_2\circ\mu_1(C_1)$ form even permutations, and the other two form odd permutations.  Of the ``even" terms, one has coefficient $-1$ and one has coefficient $+1$.  The same holds for the ``odd" terms.  Therefore, when we apply $(1+\omega_1+\omega_1^2)$ to the terms with coefficient $+1$, we get six terms with coefficient $+1$, one for each element of $S_3$.  When we apply $(1+\omega_1+\omega_1^2)$ to the terms with coefficient $-1$, we get the same six terms with negative coefficients, so all terms cancel.  To verify the above claims, one can examine all possible Gauss diagrams of $C_1$ with two non-crossing arrows, corresponding to the self-intersection points $p$ and $q$ (there are three such diagrams).   Figure 16 lists the terms in the sum $(1+\omega_i +\omega_i^2)\circ \mu_{i+1}\circ \mu_i$ before cancellations are made for a sample free loop.
\qed\\\\
\begin{figure} 
\centering
\subfigure 
{
    \label{fig:sub:a}
    \includegraphics[width=2cm]{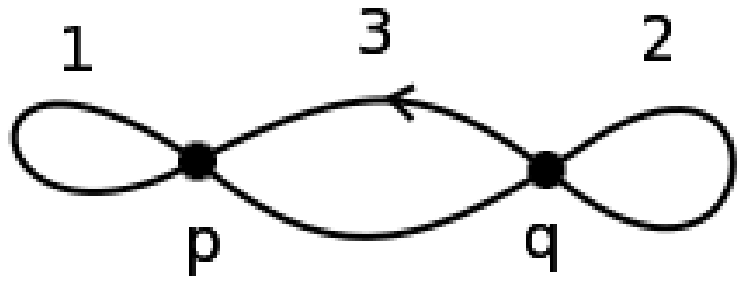}
}
\hspace{.1cm}{\bf -}
\subfigure 
{
    \label{fig:sub:a}
    \includegraphics[width=2cm]{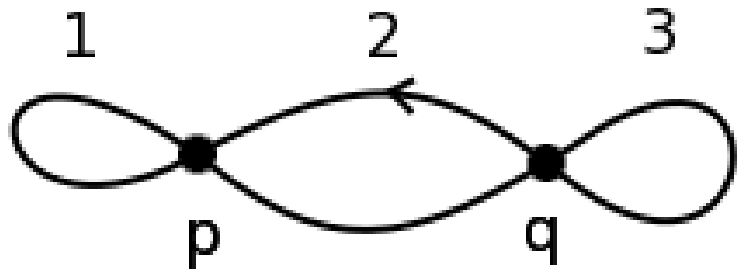}
}
\hspace{.1cm}{\bf +}
\subfigure 
{
    \label{fig:sub:a}
    \includegraphics[width=2cm]{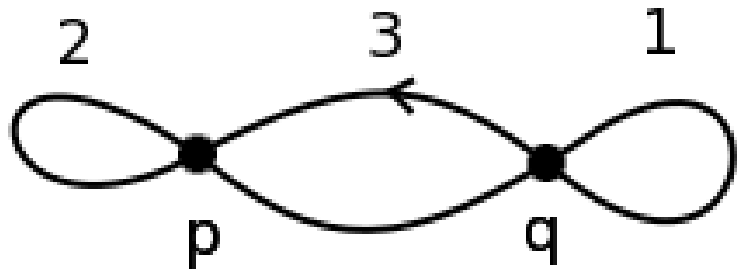}
}
\hspace{.1cm}{\bf -}
\subfigure 
{
    \label{fig:sub:b}
    \includegraphics[width=2cm]{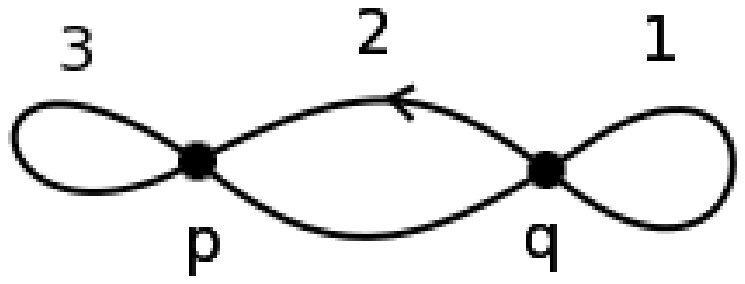}
}

\label{fig:sub} 
\end{figure}
\begin{figure}
\centering
\subfigure 
{
    \label{fig:sub:a}
    \includegraphics[width=2cm]{loop2intptsdotpandq321.eps}
}
\hspace{.1cm}{\bf -}
\subfigure 
{
    \label{fig:sub:a}
    \includegraphics[width=2cm]{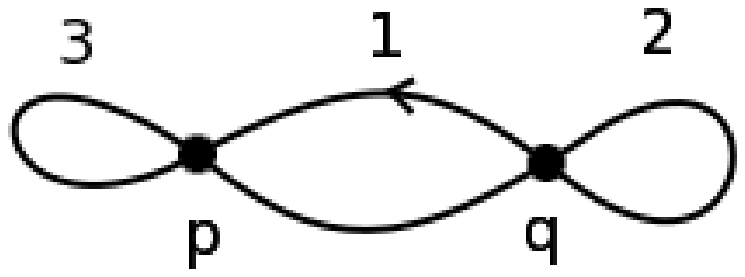}
}
\hspace{.1cm}{\bf +}
\subfigure 
{
    \label{fig:sub:a}
    \includegraphics[width=2cm]{loop2intptsdotpandq123.eps}
}
\hspace{.1cm}{\bf -}
\subfigure 
{
    \label{fig:sub:b}
    \includegraphics[width=2cm]{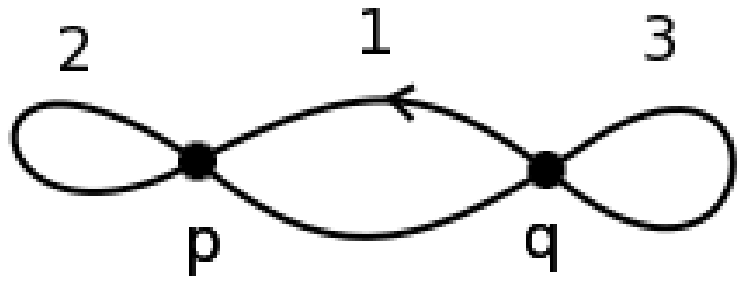}
}

\label{fig:sub} 
\end{figure}
\begin{figure}
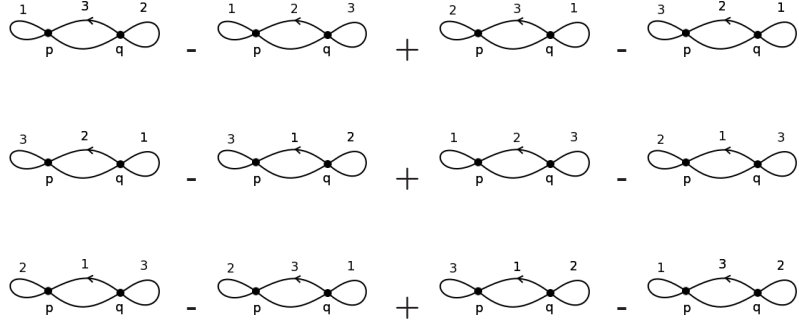

\centering
\subfigure 
{
    \label{fig:sub:a}
    \includegraphics[width=2cm]{loop2intptsdotpandq213.eps}
}
\hspace{.1cm}{\bf -}
\subfigure 
{
    \label{fig:sub:a}
    \includegraphics[width=2cm]{loop2intptsdotpandq231.eps}
}
\hspace{.1cm}{\bf +}
\subfigure 
{
    \label{fig:sub:a}
    \includegraphics[width=2cm]{loop2intptsdotpandq312.eps}
}
\hspace{.1cm}{\bf -}
\subfigure 
{
    \label{fig:sub:b}
    \includegraphics[width=2cm]{loop2intptsdotpandq132.eps}
}

\caption{The generalized co-Jacobi identity.}
\label{cojacobi.fig} 
\end{figure}

\noindent We conclude by stating the relationship between $\mu_i$ and the Turaev cobracket.  Let $E$ be a map which erases all chords from a diagram, and tensors the resulting loops, putting $C_i$ in the $i^{th}$ position of the tensor product.  Let $\Delta_i=Id\otimes ...\otimes Id \otimes \Delta \otimes Id \otimes...\otimes Id$, where $\Delta$ is in the $i^{th}$ position.  
\begin{prop}  $E\circ \mu_i=\Delta_i \circ E$. \label{smoothingprop}
\end{prop}
\noindent\textit{Proof.} This follows from the fact that $E\circ \mu_1=\Delta([\alpha])$ for any free loop $\alpha$ on $F$.
\qed\\\\
\textit{Remark:}  One might hope to find analogues of the compatibility and involutivity conditions in the Goldman-Turaev Lie bialgebra:
\begin{itemize}
\item $\Delta$ and $[,]$ satisfy the compatibility condition $\Delta\circ [\alpha,\beta]=[\alpha,\Delta(\beta)]+[\Delta(\alpha),\beta]$, where $[\alpha,\beta\otimes\gamma]$ is given by $[\alpha,\beta]\otimes \gamma+ \beta\otimes[\alpha,\gamma]$.
\item The Lie bialgebra formed by $[,]$ and $\Delta$ is involutive, i.e., $[,]\circ\Delta=0$.
\end{itemize}
However, we do not see a way of doing this unless we allow internal chords in our chord diagrams, and once we do this, it is unclear whether the maps involved in analogues of these identities (and $\mu$ in particular) are well-defined.\\\\
\noindent \textbf{Acknowledgements.}
I would like to thank my advisor Vladimir Chernov for his guidance and for reading and commenting on many drafts of this paper.  I would also like to thank Yong Su for computing examples using a preliminary definition of the operation $\mu$.

\end{document}